\documentclass[11pt, a4paper]{article}

\usepackage{bbm}
\usepackage{graphicx}
\usepackage{amsmath}
\usepackage{amsfonts}
\usepackage{amssymb}
\usepackage{enumerate}
\usepackage{lscape}
\usepackage{longtable}
\usepackage{rotating}
\usepackage{color}
\usepackage{url}
\usepackage{subfigure}
\usepackage{rotating}
\newtheorem{theorem}{Theorem}[section]

\newtheorem{algorithm}{Algorithm}[section]

\newtheorem{definition}{Definition}[section]

\newtheorem{lemma}{Lemma}[section]

\newtheorem{remark}{Remark}[section]

\newtheorem{assumption}{Assumption}[section]

\newcommand{\proof}[1]{{\noindent\emph{Proof.}} #1 \hfill$\Box$\linebreak}
\newcommand{\beqn}[1]{\begin{equation}\label{#1}}
\newcommand{\eeqn}{\end{equation}}
\def\RR{\mathbb{R}}

\def\EE{\mathbb{E}}

\newcommand{\half}{{\frac{1}{2}}}
\newcommand{\sumn}{{\sum_{k=0}^{N_\epsilon-1} }}

\newcommand{\Alpha}{{\cal A}}

\begin{document}
\title{
Global convergence rate analysis of unconstrained optimization methods based on probabilistic models }
\author{
C. Cartis\thanks{Mathematical Institute, University of Oxford, Radcliffe Observatory Quarter, Woodstock Road,
Oxford, OX2 6GG, United Kingdom ({\tt cartis@maths.ox.ac.uk}). This work was partially supported by the Oxford University EPSRC Platform Grant
 EP/I01893X/1.}
\and
K. Scheinberg\thanks{
Department of Industrial and Systems Engineering, Lehigh University,
Harold S. Mohler Laboratory, 200 West Packer Avenue, Bethlehem, PA 18015-1582, USA
({\tt katyas@lehigh.edu}). The work of this author is partially supported by NSF Grants DMS 10-16571, DMS 13-19356, CCF-1320137,
 AFOSR Grant FA9550-11-1-0239, and  DARPA grant FA 9550-12-1-0406 negotiated by AFOSR.
}
}
\maketitle
\footnotesep=0.4cm
{\small
\begin{abstract}

We present global 
convergence rates for a line-search method which is based on random
first-order models and directions whose 
quality is ensured only with certain probability. We show that in terms of the order of the accuracy, the evaluation 
complexity of such a method is the same as  its counterparts that use
deterministic accurate models; the use of 
probabilistic models only increases the complexity by a constant, which depends on the probability of the models being good. 
We particularize and improve  these results in the convex and strongly convex case. 

We also analyze a probabilistic cubic regularization 
variant that allows approximate probabilistic second-order models and show improved complexity bounds 
compared to probabilistic first-order methods; again, as a function of the accuracy, the probabilistic cubic regularization
bounds are of the same (optimal) order as for the deterministic case.



\end{abstract}

\bigskip

\begin{center}
\end{center}
}

\section{Introduction}

We consider in this paper the unconstrained optimization problem
\[
\min_{x \in \mathbb{R}^n} f(x),
\]
where the first (and second, when specified) derivatives of the objective function $f(x)$ are
assumed to exist and be (globally) Lipschitz continuous.

Most unconstrained  optimization methods rely on approximate local information to compute a local descent step 
 in such a way that sufficient decrease of the objective function is achieved. To ensure such sufficient decrease, 
 the step has to satisfy certain requirements. 
Often in practical applications ensuring these requirements for each  step is prohibitively 
expensive or impossible. This may be due to the fact that derivative information about the objective 
function is not available or because full gradient (and Hessian)  are too expensive to compute, or a 
model of the  objective function is too expensive to optimize accurately.

Recently, there has been a significant increase in interest in
unconstrained optimization methods with inexact information. Some of
these methods consider the case when gradient information is
inaccurate. This error in the gradient computation may simply be
bounded in the worst case (deterministically), see, for example,
\cite{ Devolderetal2014,Schmidtetal2011}, or the error is random and
the estimated gradient  is accurate in expectation,  as in stochastic
gradient algorithms, see for example,   \cite{GhadimiLan2013,
  RobbinsMonro, SPSA, Schmidtetal13}. These methods are typically
applied in a  convex setting and do not extend to nonconex
cases. Complexity bounds are derived that bound the expected accuracy
that is  achieved after a given number of iterations.   

In the nonlinear optimization setting, the complexity of various
unconstrained methods has been derived under exact derivative
information \cite{cagotoPI,cagotoPII,NesterovPolyak2006}, and
also under inexact information, where the errors are bounded in a
deterministic fashion \cite{Byrdetal2013, cgt39, Devolderetal2014,
  Saundersetal, Schmidtetal2011}. In all the cases of the deterministic inexact setting, traditional optimization algorithms such as line search, trust region or adaptive regularization algorithms are applied with little modification and work in practice as well as in theory, while the error is assumed to be bounded in some decaying manner {\em at each iteration}. In contrast, the methods based on stochastic  estimates of the derivatives, do not assume deterministically bounded errors,  however they are quite different from the "traditional" methods  in their strategy for step size selection and averaging of the iterates.  In other words, they are not simple counterparts of the deterministic methods. 

Our purpose in this paper is to derive a class of methods which
inherit the best properties of traditional deterministic algorithms,
and yet relax the assumption that the
derivative/model error is bounded in a deterministic manner. Moreover,
we do not assume that the error is zero in expectation or that it has
a bounded variance. Our results apply in the setting where at each
iteration, with sufficiently high probability, the error is bounded in
a decaying manner, while in the remaining cases, this error can be
arbitrarily large. In this paper, we  assume that the error may happen
in the computation of the derivatives and search directions, but that
there is no error in the function evaluations, when success of an iterate has to be validated. 

Recently several methods for unconstrained black-box optimization have been proposed, which rely on random models or directions \cite{Bandeiraetal2014,Grattonetal2014, YNesterov_2011}, but are applied to deterministic functions. In this paper we take this line of work one step further by establishing expected convergence rates for several schemes based on one generic analytical framework. 

We consider four cases and derive four different complexity bounds. In
particular, we analyze a line search method based on random models,
for the cases of general nonconex, convex and strongly convex
functions. We also analyze a second order method - an adaptive
regularization method with cubics \cite{cagotoPI, cagotoPII} - which
is known to achieve the optimal convergence rate for the nonconvex smooth
functions \cite{cgt38} and we show that the same convergence rate holds in
expectation. 

In summary, our results differ from existing literature using inexact, stochastic or random information in the following main points:
\begin{itemize}
\item Our models are assumed to be "good" with some probability, but there is no other assumptions on the expected values or variance of the model parameters. 

\item The methods that we analyze are essentially the exact counterparts of the deterministic methods, and do not require averaging of the iterates or any other significant changes. We believe that, amongst other things, our analysis helps to understand the convergence properties of practical algorithms, that do not always seek to ensure theoretically required model quality. 

\item Our main convergence rate results provide a bound on the {\it
    expected number of iterations } that the algorithms take before
  they achieve a desired level of accuracy. This is in contrast to a
  typical analysis of randomized or stochastic methods, where what is
  bounded is the expected error after a given number of
  iterations. Both bounds are useful, but we believe that the bound on
  the expected number of steps is a somewhat more meaningful
  complexity bound in our setting. 
 The only other  work that we are aware of which provides bounds in terms of the number of required steps is \cite{Grattonetal2014} where 
probabilistic bounds are derived in the particular context of random direct search with possible extension to trust region methods  as discussed in Section 6 of \cite{Grattonetal2014}.

\end{itemize}

An additional  goal of this paper is to present a general theoretical
framework, which could be used to analyze the behavior of other algorithms,
and different possible model construction mechanisms under the assumption that the objective function is deterministic. 
We propose a general analysis of an optimization scheme by reducing it to the analysis of a stochastic process. 
Convergence results for a trust region method in \cite{Bandeiraetal2014} also rely on a stochastic process analysis, but 
only in terms of  behavior in the limit. These results have now been extended to noisy (stochastic) functions, see \cite{RChen_2015, Chenetal2015}.  
Deriving convergence {\em rates} for methods applied to  stochastic functions is the subject of future work and is likely to depend on the results in this paper.

The rest of the paper is organized as follows. In Section \ref{sec:genralscheme}
we describe the general scheme which encompasses several unconstrained optimization methods. This scheme is based on using random models, which are assumed to satisfy some "quality" conditions with probability at least $p$, conditioned on the past. Applying this optimization scheme results in a stochastic process, whose behavior is analyzed in the later parts of Section \ref{sec:genralscheme}. Analysis of the stochastic process allows us to bound the expected number of steps of our generic scheme until a desired accuracy is reached. 
In Section \ref{sec:ls_algorithm} we analyze a linesearch algorithm
based on random models and show how its behavior fits into our general
framework for the cases of nonconex, convex and strongly convex
functions. In Section \ref{sec:arc_algorithm} we apply our generic
analysis to the case of the Adaptive Regularization method with Cubics
(ARC). Finally, in Section \ref{sec:models} we describe different 
settings where the models of the objective functions satisfy the probabilistic conditions of our schemes. 


\section{A general optimization scheme with random models}\label{sec:genralscheme}

This section presents the main features of our algorithms and analysis,
in a general framework that we will, in subsequent sections,
particularize to specific algorithms (such as linesearch and cubic
regularization) and classes of functions (convex, nonconvex). The
reasons for the initial generic approach is to avoid repetition of
the common elements of the analysis for the different algorithms and to emphasize the key ingredients
of our analysis, which is possibly applicable to other algorithms
(provided they satisfy our framework).

\subsection{A general optimization scheme}

We first describe a generic algorithmic framework that encompasses the
main components of the unconstrained optimization schemes we analyze
in this paper. The scheme relies on building a model of the objective
function at each iteration, 
minimizing this model or reducing it in a sufficient manner and
considering the step which is dependent on a stepsize parameter and
which provides the model reduction (the stepsize parameter may be
present in the model or independent of it). 
This step determines a new candidate point. The function value is then
computed (accurately) at the new candidate point. 
If the function reduction provided by the candidate point is deemed
sufficient, then the iteration is declared successful, the 
candidate point becomes the new iterate and the step size parameter is
increased. Otherwise, the iteration is 
unsuccessful, the iterate is not updated and the step size parameter is reduced.

We summarize the main steps of the scheme below.

\begin{algorithm}\label{alg:generic}{\bf Generic optimization framework based on  random models} 
 \ \\
\begin{description}
\item[Initialization] \ \\
Choose a class of  (possibly random) models $m_k(x)$, choose constants $\gamma\in (0,1)$,
$\theta\in (0,1)$, $\alpha_{\max}>0$.
Initialize the algorithm by choosing $x_0$, $m_0(x)$, $0<\alpha_0<\alpha_{\max}$. 

 \item[1. Compute  a model and a step] \ \\
Compute a local (possibly random) model $m_k(x)$ of $f$ around $x^k$.\\
Compute a step $s^k(\alpha_k)$ which reduces $m_k(x)$, where the
parameter $\alpha_k>0$ is present in the model or  in the
step calculation. 

\item[2. Check sufficient decrease]\ \\  
Compute $f(x^k+s^k(\alpha_k))$ 
and check if sufficient reduction (parametrized by $\theta$) is achieved in $f$ with respect to
${m_k(x^k)-m_k(x^k+s^k(\alpha_k))}$.

\item[3. Successful step] \ \\ 
If sufficient reduction is achieved then, $x^{k+1}:=x^k+s^k(\alpha_k)$,  set $\alpha_{k+1}= 
\min\{\alpha_{\max},\gamma^{-1}\alpha_k\}$. Let $k:=k+1$.

\item[4. Unsuccessful step]\ \\  
Otherwise, $x^{k+1}:=x^k$, set $\alpha_{k+1}=\gamma \alpha_k$. Let $k:=k+1$.\\
 \end{description}
 \end{algorithm}

Let us illustrate how the above scheme relates to standard
optimization methods. In {\it linesearch methods}, one minimizes a
linear model $m_k(x)=f(x^k)+  (x-x^k)^Tg^k$ (subject to some normalization), or a quadratic one
$m_k(x)=f(x^k)+  (x-x^k)^Tg^k+\frac{1}{2}(x-x^k)^\top b^k(x-x^k)$
(when the latter is well-defined, with $b^k$ - a Hessian approximation matrix), 
to find directions $d^k=-g^k$ or $d^k=-(b^k)^{-1}g^k$,
respectively. Then the step is defined as $s^k(\alpha_k)=\alpha_kd^k$ for some
$\alpha_k$ and, commonly,   the (Armijo) decrease condition is checked,
\[
 f(x^k)-f(x^k+s^k(\alpha_k))\geq - \theta s^k(\alpha_k)^Tg^k,
 \]
where $- \theta s^k(\alpha_k)^Tg^k$ is a multiple of 
$m_k(x^k)-m_k(x^k+s^k(\alpha_k))$.  Note that if  the model stays the same in that $m_k(x) \equiv
m_{k-1}(x)$ for each $k$, such that  $(k-1)$st iteration is unsuccessful, then 
the above framework essentially reduces to a standard deterministic linesearch. 

 In the case of {\it cubic regularization methods}, $s^k(\alpha_k)$ is
 computed to approximately minimize a cubic model $m_k(x)=f(x^k)+
 (x-x^k)^Tg^k+\frac{1}{2}(x-x^k)^\top
 b^k(x-x^k)+\frac{1}{3\alpha_k}\|x-x^k\|^3$ and the sufficient decrease condition is 
 \[
\frac{f(x^k)-f(x^k+s^k(\alpha_k))}{m(x^k)-m(x^k+s^k(\alpha_k))} \geq \theta>0. 
\]
Note that here as well, in the deterministic case, $g^k=g^{k-1}$ and
$b^k=b^{k-1}$ for each $k$ such that $(k-1)$st iteration is unsuccessful but
$\alpha_k\neq \alpha_{k-1}$.

The key assumption in the usual deterministic case is that the models $m_k(x)$ are
sufficiently accurate in a small neighborhood of the current iterate
$x^k$. The goal of this paper is to relax this requirement and allow
the use of random local models which are accurate only with certain
probability (conditioned on the past). In that case, note that the models need to be re-drawn
after each iteration, whether successful or not. 

Note that our general setting includes the cases when the model (the derivative information, for example) is always accurate, but the step $s^k$ is computed approximately, in a probabilistic manner. For example,
$s^k$ can be an approximation of $-(b^k)^{-1}g^k$. It is easy to see how randomness in $s^k$ calculation can be viewed as the randomness in the model, by considering that instead of the accurate model
\[
f(x^k)+  (x-x^k)^Tg^k+\frac{1}{2}(x-x^k)^\top b^k(x-x^k)
\]
we use an approximate model
\[
m_k(x)=f(x^k)-  (x-x^k)^Tb^ks^k+\frac{1}{2}(x-x^k)^\top b^k(x-x^k).
\]
Hence, as long as the accuracy requirements are carried over accordingly the approximate random models subsume the case of approximate random step computations. The next section
makes precise our requirements on the probabilistic models.

\subsection{Generic probabilistic models}

We will now introduce the key probabilistic ingredients of our scheme. In particular we assume that our models $m_k$ are random and that they satisfy some notion of good quality 
with some probability $p$.  We will consider random models $M_k$, and then use the notation
$m_k= M_k(\omega_k)$ for their realizations. The randomness of the models will imply the randomness  of the points~$x^k$, the step length parameter~$\alpha_k$, the computed steps $s^k$ and  other quantities produced by the algorithm. 
Thus, in our paper, these random variables will be denoted by $X^k$, $\Alpha_k$, $S^k$ and so on, respectively,
while $x^k= X^k(\omega_k)$, $\alpha_k= \Alpha_k(\omega_k)$,  $s^k= S^k(\omega_k)$, etc, denote their realizations (we will omit the $\omega_k$ in the notation for brevity). 

For each specific optimization method, we will define a notion of
sufficiently accurate models. The desired accuracy of the model
depends  on the current iterate $x^k$, step 
parameter $\alpha_k$ and,
possibly, the step $s^k(\alpha_k)$. This notion involves model
properties which make sufficient decrease in $f$ 
achievable by the step $s^k(\alpha_k)$. Specific conditions on the
models will be stated for each algorithm in the 
respective sections and how these conditions may be achieved will be discussed in Section \ref{sec:models}. 

\begin{definition}\label{def:random_poised_mart} {\bf [sufficiently
    accurate models; true and false iterations]}
We say that a sequence of random models $\{M_k\}$ is $(p)$-probabilistically
 ``sufficiently accurate" for a corresponding sequence $\{\Alpha_k, X^k\}$, if the following indicator random variable
\[
I_k \; = \; \mathbbm{1}\{M_k \text{ is a sufficiently accurate  model of }f \text{ for the given  }X^k {\text and\ } \Alpha_k\}
\]
satisfy the following submartingale-like condition
\begin{equation}\label{prob-p-def}
P(I_k=1|\,  {\cal F}^M_{k-1}) \; \geq \; p,
\end{equation}
where ${\cal F}^M_{k-1}=\sigma(M_0,\ldots,M_{k-1})$ is the
$\sigma$-algebra generated by $M_0,\ldots,M_{k-1}$ - in other words,
the history of the algorithm up to iteration $k$. 

We say that iteration $k$ is a {\bf true} iteration if the event $I_k=1$ occurs. Otherwise the iteration is called {\bf false}.
\end{definition}

Note that $M_k$ is a random model that, given the past history, encompasses all the randomness of iteration $k$ of our algorithm. The iterates $X^k$ and the step length parameter$\Alpha_k$ are random
variables defined over the $\sigma$-algebra generated by $M_0,\ldots,M_{k-1}$.  Each $M_k$ depends on $X^k$ and $\Alpha_k$ and hence on $M_0,\ldots, M_{k-1}$.
Definition~\ref{def:random_poised_mart} serves to enforce the following property:  even though the accuracy of $M_k$ may be dependent on the history, ($M_1,\ldots,M_{k-1}$), via its dependence on $X^k$ and $\Alpha_k$, it is sufficiently good  with probability at least~$p$, regardless of that history. This condition is more reasonable  than   complete independence of $M^k$ from the past, which is difficult to ensure. 
It is important to note that, from this assumption, it follows that whether or not the step is deemed successful and the iterate $x^k$ is updated, our scheme always updates the model $m_k$, unless $m_k$ is somehow known to be sufficiently accurate for $x^{k+1}=x^k$ and $\alpha_{k+1}$. We will discuss this in more detail in Section \ref{sec:models}. 

When Algorithm \ref{alg:generic} is  based on probabilistic models
(and all its specific variants under consideration), it results in a
discrete time stochastic process. This stochastic process encompasses
random elements such $\Alpha_k$, $X^k$, $S^k$, which are directly computed by the algorithm, but also some quantities that can be derived as functions of $\Alpha_k$, $X^k$, $S^k$, such as $f(X^k)$, $\|\nabla f(X^k)\|$ and a quantity $F_k$, which we will use to denote some measure of progress towards optimality.  
Each realization of the sequence of random models results in a
realization of the algorithm, which in turn produces the corresponding
sequences $\{\alpha_k\}$, $\{x^k\}$, $\{s^k\}$,  $\{f(x^k)\}$,
$\{\|\nabla f(x^k)\|\}$ and $\{f_k\}$\footnote{Note that throughout,
  $f(x^k)\neq f_k$, since $f_k$ is a related measure of progress towards optimality.}. 
  We will analyze the stochastic  processes 
  restricting our attention to some of the random quantities that belong to this process and will ignore the rest, for the brevity of the presentation. Hence when we say that Algorithm 
\ref{alg:generic} generates the stochastic process $\{X^k, \Alpha_k\}$, this means we want to focus on the properties of these random variables, but keeping in mind that there are other random quantities in this stochastic process.

We will derive complexity bounds for each algorithm in the following
sense. We will define the accuracy goal that we aim to reach and then
we will bound the expected number of steps that the algorithm takes
until this goal is achieved.  The analyses will
follow common steps, and the main ingredients are described below.
We then apply these steps to each case under consideration. 

\subsection{Elements of global convergence rate analysis} \label{sec:elem_conv_anal}

First we recall a standard notion from stochastic processes. 
\paragraph*{Hitting time.}
For a given discrete time stochastic process, $Z_t$, recall the concept of a
{\em hitting time} for an event $\{Z_t\in S\}$. This is a random variable, defined as $T_S=\min \{t:\
Z_t\in S\}$ - the first time the event  $\{Z_t\in S\}$ occurs. In our
context, set $S$ will either be a set of real numbers larger than some
given value, or smaller than some other given value. 

\paragraph{Number of iterations $N_{\epsilon}$ to reach $\epsilon$ accuracy.}
\label{Neps-page}
Given a level of accuracy $\epsilon$, we aim to derive a bound on the expected number of iterations $\EE(N_\epsilon)$ which occur in the algorithm until the given accuracy level is reached.
The number of iterations $N_{\epsilon}$ is a random variable, which
can be defined as a hitting time of some stochastic
   process, 
dependent on the case under analysis.
In particular, 

\begin{itemize}
\item If $f(x)$ is not known to be convex, then
  $N_\epsilon$ is  the hitting time
   for $\{\|\nabla f(X_k)\|\leq \epsilon\}$, namely,
the number of steps  the algorithm takes until  $\|\nabla f (X^k)\|\leq \epsilon$ occurs for the first time. 
 
\item
If $f(x)$ is convex or strongly convex then
$N_\epsilon$ is the hitting time for
$\{f(X^k)-f_*\leq \epsilon\}$, namely,
 the number of steps  the algorithm takes until  $f(X^k)-f_*\leq \epsilon$ occurs for the first time, where $f_*=f(x^*)$ with $x^*$, a global minimizer of $f$. 
\end{itemize}


We will bound $\EE(N_\epsilon)$ by  observing that for all $k<N_\epsilon$ the 
stochastic process induced by Algorithm \ref{alg:generic} behaves in a certain way. 
 To formalize this, we need to define the following random variable and its upper
bound.  

\paragraph{Measure of progress towards optimality,  $F_k$.} This  measure
is  defined by the  total function decrease or by the distance to the
optimum. In particular,
\begin{itemize}
\item   If $f(x)$ is not known to be convex, then
  $F_k=f(X^0)-f(X^k)$. 
\item
If  $f(x)$ is convex, then $F_k=1/(f(X^k)-f_*)$. 
\item
If  $f(x)$ is strongly convex, then $F_k=\log(1/(f(X^k)-f_*))$. 
\end{itemize}

\paragraph{Upper bound  $F_{\epsilon}$ on $F_k$.}
From the algorithm construction,  $F_k$  defined above is always
nondecreasing and there exists a deterministic upper bound $F_{\epsilon}$ 
in each case, defined as follows.
\begin{itemize}
\item   If $f(x)$ is not known to be convex, then
  $F_\epsilon=f(X^0)-f_*$, where $f_*$ is a global lower bound on $f$.
\item
If  $f(x)$ is convex, then $F_\epsilon=1/ \epsilon$. 
\item
If  $f(x)$ is strongly convex, then $F_\epsilon=\log(1/ \epsilon)$. 
\end{itemize}

We observe that $F_k$ is a nondecreasing process and $F_\epsilon$ is the largest
possible value  that $F_k$ can achieve. 

Our analysis will  be based on the following observations, which
are borrowed from the global rate  analysis of the deterministic
methods \cite{nesterovtext}. 

\begin{itemize}
\item {\bf Guaranteed amount of increase in $f_k$.}  For all $k<N_\epsilon$ (i.e., until  
the desired accuracy has  been reached), if the $k$th iteration
  is true and successful, 
 then $f_k$ is increased by an amount proportional to $\alpha_k$.

\item {\bf Guaranteed threshold for $\alpha_k$.} There exists a
  constant, which we will call $C$, such that if $\alpha_k  \leq  C$ 
and the $k$th iteration is true, 
 then the $k$th iteration is also successful,  and hence 
$\alpha_{k+1}=\gamma^{-1}\alpha_k$. This constant $C$ depends on the
algorithm and Lipschitz constants of $f$. 

\item {\bf Bound on the number of  iterations.}
If all iterations were true, then by the above observations, 
$\alpha_{k}\geq \gamma C$ and, hence, $f_k$ increases by at
least a constant for all $k$. From this a bound on the number of
iterations, knowing that $f^k$ cannot exceed $F_\epsilon$. 
\end{itemize}
 In our case
not all iterations are true, however, under the assumption that they ``tend'' to be true,
 as we will show, when
$\Alpha_{k}\leq C$, then iterations ``tend'' to be successful, $\Alpha_k$
``tends'' to stay near the value $C$   and the values
$F_k$ ``tend'' to increase by a constant. The analysis is then performed
via a study of stochastic processes, which we describe in detail
next.

\subsection{Analysis of the  stochastic processes }\label{sec:akfkproc}
Let us  consider the stochastic process $\{\Alpha_k, F_k\}$ generated by Algorithm \ref{alg:generic} using random, 
$p$-probabilistically sufficiently accurate models $M_k$, with $F_k$ defined above. 
Under the assumption that the sequence of models $M_k$ are $p$-probabilistically sufficiently accurate, 
each iteration is true with probability at least $p$, conditioned on the past. 

We assume now (and we show later for each specific case)  that $\{\Alpha_k, F_k\}$ obeys the following rules for all $k<N_\epsilon$.
\begin{assumption}\label{ass:alg_behave}
There exist a constant $C>0$ and a nondecreasing function $h(\alpha)$,  $\alpha\in \mathbb{R}$,   which satisfies $h( \alpha)>0$ for any $\alpha>0$, such that for any realization of  Algorithm 
 \ref{alg:generic} the following hold for all $k<N_\epsilon$:
\begin{itemize}
\item[(i)]   If iteration $k$ is true (i.e.  $I_k=1$) and successful, then $f_{k+1}\geq f_k+h(\alpha_k)$. 

\item[(ii)] If $\alpha_k  \leq  C$ and iteration $k$ is true  then
  iteration $k$ is also successful, which implies
$\alpha_{k+1}=\gamma^{-1}\alpha_k$.
\item[(iii)] $f_{k+1}\geq f_k$ for all $k$.
\end{itemize}
\end{assumption}

For future use let us state an auxiliary lemma.

\begin{lemma}\label{lem:frac_of_true}
Let $N_\epsilon$ be the hitting time  as defined on page \pageref{Neps-page}. For all $k<N_\epsilon$,
 let $I_k$ be the sequence of random variables in Definition \ref{def:random_poised_mart} so that \eqref{prob-p-def} holds.
Let $W_k$ be a nonnegative stochastic process such that $\sigma(W_k)\subset {\cal F}^M_{k-1}$, 
for any $k\geq 0$. Then
\[
\EE\left (\sumn  W_kI_k \right )\geq p \EE \left (\sumn W_k\right ).
\]
Similarly,
\[
\EE\left (\sumn  W_k(1-I_k )\right )\leq (1- p) \EE \left (\sumn W_k\right ).
\]
\end{lemma}
\proof{ The proof is a simple consequence of properties of expectations, see for example, \cite[property H$^*$, page 216]{shiryaev},
\[
\EE (I_k|\, W_k) = \EE (\EE(I_k|\,  {\cal F}^M_{k-1})|\, W_k) \geq \EE (p|\, W_k) \geq p,
\]
where we also used  that $\sigma(W_k)\subset {\cal F}^M_{k-1}$. 
Hence by the law of total expectation, we have $\EE(W_kI_k)=\EE(W_k\EE(I_k|W_k))\geq p\EE(W_k)$. Similarly, 
we can derive $\EE(\mathbbm{1}\{k<N_\epsilon\}W_kI_k)\geq p\EE(\mathbbm{1}\{k<N_\epsilon\}W_k)$, because $\mathbbm{1}\{k<N_\epsilon\}$ is also determined
by  ${\cal F}^M_{k-1}$. Finally,
\[
\EE\left (\sumn W_kI_k\right )= \EE\left (\sum_{k=0}^\infty \mathbbm{1}\{k<N_\epsilon\}W_kI_k\right )\geq p\EE\left (\sum_{k=0}^\infty \mathbbm{1}\{k<N_\epsilon\}W_k\right )=
p\EE\left (\sumn W_k\right ).
\]
The second inequality is proved analogously. 
}

Let us now define two indicator random variables, in addition to $I_k$ defined earlier,
\[
\Lambda_k=\mathbbm{1}\{ \Alpha_k >C\}, 
\]
and 
\[
\Theta_k=\mathbbm{1}\{{\rm Iteration\ } k\ {\rm is\ successful\ i.e.,\ } \Alpha_{k+1}=\gamma^{-1}\Alpha_k \}. 
\]
Note that $\sigma(\Lambda_k)\subset {\cal F}_{k-1}^M$ and $\sigma(\Theta_k)\subset {\cal F}_{k}^M$, that is the random variable 
$\Lambda_k$ is fully determined by the first $k-1$ steps of the algorithm, while $\Theta_k$ is fully determined by the first $k$ steps. 
We will use $\lambda_k$, $i_k$ and $\theta_k$ to denote realizations of $\Lambda_k$, $I_k$ and $\Theta_k$, respectively. 

These indicators will help us define our algorithm more  rigorously as a stochastic process. 
Without loss of generality, we assume that $C=\gamma^{c}\alpha_0<\gamma\alpha_{\max}$ for some positive integer $c$. In other words, $C$ is the largest value that the step size $\Alpha_k$ actually achieves for which  part $(ii)$ of Assumption \ref{ass:alg_behave} holds.
The condition $C< \gamma\alpha_{\max}$ is a simple technical condition,
which is not necessary, but which simplifies the presentation later in this section. 
Under Assumption \ref{ass:alg_behave}, recalling the update rules for
$\alpha_k$ in  Algorithm \ref{alg:generic}  and the assumption that
true iterations occur with probability
at least $p$,   we can write the stochastic process $\{\Alpha_k, F_k\}$ as obeying the expressions below:

\begin{equation}\label{eq:proc1_Zk}
\Alpha_{k+1}=\left \{ \begin{array}{ll} \gamma^{-1} \Alpha_k& 
{\rm if\  } I_k=1\ {\rm and \ } \Lambda_k=0,
\\ \gamma \Alpha_k & {\rm if\  }I_k=0\ {\rm and \ } \Lambda_k=0, \\
\min\{\alpha_{\max},\gamma^{-1} \Alpha_k \}& {\rm if\  }  \Theta_k=1\ {\rm and \ } \Lambda_k=1, \\
 \gamma \Alpha_k &    {\rm if\  }  \Theta_k=0\ {\rm and \ } \Lambda_k=1,
\end{array}\right . 
\end{equation}
 \begin{equation}\label{eq:proc1_Yk}
F_{k+1}\geq \left \{ \begin{array}{ll} F_k+h(\Alpha_k) & {\rm if\  } I_k=1\ {\rm and \ } \Lambda_k=0,
 \\ F_k & {\rm if\  }I_k=0\ {\rm and \ } \Lambda_k=0,\\
F_k+h(\Alpha_k)
&   \Theta_kI_k=1 \ {\rm and \ } \Lambda_k=1, \\ F_k &   \Theta_k I_k=0 \ {\rm and \ } \Lambda_k=1. \end{array}\right . 
\end{equation}

%

We conclude that, when  $\Alpha_k\leq C$, 
a  successful iteration happens with probability at least $p$, and in
that case $\Alpha_{k+1}=\gamma^{-1}\Alpha_k$, and that an unsuccessful iteration happens with probability at most
  $1-p$, in which case $\Alpha_{k+1}=\gamma\Alpha_k$. 
  Note that there is no known probability bound for the different
  outcomes when $\Alpha_k>C$. However,  we know that $I_k=1$ with probability at least $p$ and if, in addition,
    iteration $k$ happens to be successful, then $F_k$ is increased  by at least  $h(\Alpha_k)$. 

In summary, from the above discussion, we have 
\[
\begin{array}{l}{\textit{for all $k<N_{\epsilon}$, Algorithm \ref{alg:generic} under Assumption \ref{ass:alg_behave} yields}}\\
{\textit{the stochastic process 
$\{\Alpha_k, F_k\}$ in \eqref{eq:proc1_Zk} and \eqref{eq:proc1_Yk}.}}
\end{array}
\]

%
%

\subsection{Bounding the number of steps for which $\alpha_k\leq C$}\label{sec:zkykprime}

In this subsection we derive a bound on $\EE\left (\sumn (1-\Lambda_k)\right )$. The bound for $\EE(\sumn \Lambda_k)$ will be derived in the next section. 

The following simple result holds for every realization of the algorithm and  stochastic process $\{\Lambda_k, I_k, \Theta_k\}$.

\begin{lemma}\label{lem:bound_on_small}
For any $l\in \{0,\ldots, N_{\epsilon}-1\}$ and for all realizations of Algorithm \ref{alg:generic}, we have 
\[
\sum_{k=0}^l (1-\Lambda_k) \Theta_k  \leq \half (l+1).
\]
\end{lemma}
\proof{
By the definition of $\Lambda_k$ and $\Theta_k$ we know that when $(1-\Lambda_k) \Theta_k=1$ then we have a successful iteration and $\Alpha_k\leq C$.
In this case $\Alpha_{k+1}=\gamma^{-1}\Alpha_k$. It follows that amongst all iterations, at most half can be successful and have $\Alpha_k\leq C$,  because for each such iteration, when $\Alpha_k$ gets increased by a factor of $\gamma^{-1}$, there has to be at least one iteration when  $\Alpha_k$ is decreased by the same factor, since $\Alpha_0\geq C$.
}

Using this we derive the bound. 
 \begin{lemma}\label{lem:hittime1}
 \[
  \EE\left (\sumn (1-\Lambda_k)\right)\leq \frac{1}{2p}\EE(N_\epsilon)
  \]
  \end{lemma}
  \proof{
By Lemma \ref{lem:frac_of_true} applied to $W_k=1-\Lambda_k$ we have
\begin{equation}\label{eq:bound_small_1}
\EE\left (\sumn (1-\Lambda_k)I_k\right)\geq p \EE\left (\sumn (1-\Lambda_k)\right).
\end{equation}
From the fact that all true iterations are successful when $\alpha_k\leq C$, 
\begin{equation}\label{eq:bound_small_2}
\sumn (1-\Lambda_k)I_k\leq \sumn (1-\Lambda_k)\Theta_k.
\end{equation}
Finally, from Lemma \ref{lem:bound_on_small}
\begin{equation}\label{eq:bound_small_3}
\sumn (1-\Lambda_k)I_k\leq \half N_\epsilon. 
\end{equation}

Taking expectations in \eqref{eq:bound_small_2} and \eqref{eq:bound_small_3} and combining with \eqref{eq:bound_small_1}, we obtain
the result of the lemma.
}

\subsection{Bounding the expected number of steps for which $\alpha_k>C$}\label{sec:zkykprimeprime}

Let us now consider the bound on $\EE \left (\sumn \Lambda_k\right )$.  We introduce the additional notation $\bar \Lambda_k= \mathbbm{1}\{ \Alpha_k >C\}+\mathbbm{1}\{ \Alpha_k =C\}$. In other words $\bar \Lambda_k=1$ when  either $\Lambda_k=1$ or $\Alpha_k=C$.
We now define:
\begin{itemize}
\item $N_1=\sumn \bar \Lambda_k (1-I_k)\Theta_k$, which is the number of false successful iterations, when $\Alpha_k\geq C$.
\item $M_1=\sumn \bar \Lambda_k (1-I_k)$, which is the number of false  iterations, when $\Alpha_k\geq C$.
\item $N_2=\sumn \bar \Lambda_k I_k \Theta_k$, which is the number of true successful iterations, when $\Alpha_k\geq C$.
\item $M_2=\sumn \bar \Lambda_k I_k$, which is the number of true iterations, when $\Alpha_k\geq C$.
\item $N_3=\sumn \Lambda_k I_k (1-\Theta_k)$, which is the number of true unsuccessful iterations, when $\Alpha_k> C$.
\item $M_3=\sumn \Lambda_k  (1-\Theta_k)$, which is the number of unsuccessful iterations, when $\Alpha_k> C$.
\end{itemize}


Since $\EE\left (\sumn \Lambda_k\right )= \EE \left( \sumn \Lambda_k (1-I_k)\right )+\EE \left( \sumn \Lambda_k I_k\right )\leq \EE(M_1)+\EE(M_2)$, 
our goal is to bound $\EE(M_1)+\EE(M_2)$. 

Our next observation is simple but central in our analysis. It reflects the fact that the gain in $F_k$ is bounded from above by $F_\epsilon$ 
and when $\Alpha_k\geq C$ this gain is bounded from below as well, hence allowing us to bound the total number of true successful iterations when 
$\Alpha_k\geq C$.  The following two lemmas holds for every realization. 
\begin{lemma}\label{lem:bound_on_big2}
For any $l\in \{0,\ldots,N_{\epsilon}-1\}$ and for all realizations of Algorithm \ref{alg:generic}, we have 
\[
\sum_{k=0}^l  \bar \Lambda_kI_k\Theta_k \leq \frac{F_\epsilon}{h(C)},
\]
and so
\begin{equation}\label{eq:N2}
N_2\leq \frac{F_{\epsilon}}{h(C)}.
\end{equation}
\end{lemma}
\proof{
Consider any $k$ for which $\Lambda_kI_k\Theta_k =1$. From Assumption \ref{ass:alg_behave} we know that whenever an iteration is true and successful then $F_k$ get increased by at least $h(\Alpha_k)\geq h(C)$, since  $\Alpha_k\geq C$ and $h$ is nondecreasing. We also know that on other iterations $F_k$ does not decrease. The bound  $F_k\leq F_\epsilon$ trivially gives us the desired result. 
}

Another key observation is that  
\begin{equation}\label{eq:M2}
M_2\leq N_2+N_3\leq N_2+M_3,
\end{equation}
where the first inequality follows from the fact that for all $k<N_{\epsilon}$ and for all realizations, $(\bar \Lambda_k - \Lambda_k)I_k(1-\Theta_k)=0$, in other words there are no true unsuccessful iterations when $\Alpha_k=C$. 

\begin{lemma}\label{lem:bound_on_big}
For any $l\in \{0,\ldots,N_{\epsilon}-1\}$ and for all realizations of Algorithm \ref{alg:generic}, we have 
\[
\sum_{k=0}^l \Lambda_k(1- \Theta_k) \leq \sum_{k=0}^l\bar \Lambda_k\Theta_k + \log_{\gamma}\left (\frac{C}{\alpha_0}\right )
\]
%
\end{lemma}
\proof{$\Alpha_k$ is increased on successful iterations and decreased on unsuccessful ones. Hence the total number of steps
when $\Alpha_k>C$ and $\Alpha_k$ is decreased, is bounded by the total number of steps when $\Alpha_k\geq C$ is increased plus  the number of steps it is required to reduce $\Alpha_k$ from its initial value $\alpha_0$ to $C$. 
}

From Lemma \ref{lem:bound_on_big} applied to $l=N_\epsilon-1$, we can deduce that
\begin{equation}\label{eq:inc_dec}
M_3\leq N_1+N_2 + \log_\gamma(C/\alpha_0).
\end{equation}

We also have the following  lemma.
\begin{lemma}\label{no_of_false_lemma}
\begin{equation}\label{false}
\EE(M_1)\leq \frac {1-p}{p}\EE(M_2). 
\end{equation}
\end{lemma}
\proof{ 
By applying both inequalities in Lemma \ref{lem:frac_of_true} with $W_k=\bar\Lambda_k$, we obtain
\[
 \EE\left(\sumn \bar \Lambda_kI_k\right)\geq  p\EE\left(\sumn \bar\Lambda_k\right)
 \]
 and 
 \[
 \EE\left(\sumn \bar\Lambda_k(1-I_k)\right)\leq (1-p)\EE\left(\sumn \bar\Lambda_k\right)
 \]
which gives us 
\[
\EE\left (\sumn \bar\Lambda_k(I_k-1)\right ) \leq \frac {1-p}{p} \EE\left(\sumn \bar\Lambda_kI_k\right ). 
\]
}

\begin{lemma}\label{lem:hittime2}
Under the condition that $p>1/2$, we have
\[
\EE\left (\sumn \Lambda_k \right )\leq \frac{2F_{\epsilon} }{h(C) (2p-1)}+\frac{\log_\gamma(C/\alpha_0)}{2p-1}. 
\]
\end{lemma}

\proof{
Recall that $\EE\left (\sumn \Lambda_k\right )=\EE(M_1+M_2)$. 
Using \eqref{eq:M2} and \eqref{false}  it follows that
\begin{equation}\label{eq:binom}
\EE(N_1)\leq \EE(M_1)\leq\frac {1-p}{p}\EE(M_2)\leq \frac {1-p}{p}\EE(N_2+M_3)= \frac {1-p}{p}[\EE(N_2)+\EE(M_3)].
\end{equation}

Taking into account \eqref{eq:inc_dec} and using the bound  \eqref{eq:N2} on $N_2$  we have
\begin{equation}\label{eq:boundM3}
\EE(M_3)\leq \EE(N_1)+\EE(N_2) + \log_\gamma(C/\alpha_0)\leq \EE(N_1)+F_{\epsilon}/h(C) + \log_\gamma(C/\alpha_0). 
\end{equation}
Plugging this into  \eqref{eq:binom} and using the bound  \eqref{eq:N2}  on $N_2$ again, we obtain
\[
\EE(N_1)\leq  \frac {1-p}{p}\left[\frac{F_{\epsilon}}{h(C)}+\EE(N_1)+\frac{F_{\epsilon}}{h(C)}+ \log_\gamma\left(\frac{C}{\alpha_0}\right)\right],
\]
and, hence,
\[
\frac {2p-1}{p}\EE(N_1)\leq \frac {1-p}{p}\left[ \frac{2F_{\epsilon}}{h(C)} + \log_\gamma\left(\frac{C}{\alpha_0}\right)\right].
\]
This finally implies 
\begin{equation}\label{eq:boundN1}
\EE(N_1)\leq \frac {1-p}{2p-1}\left [\frac{2 F_{\epsilon}}{h(C)}+ \log_\gamma\left(\frac{C}{\alpha_0}\right)\right ].
\end{equation}
Now we can bound the expected total  number of  iterations when $\alpha_k>C$, using \eqref{eq:N2}, \eqref{eq:boundM3} and \eqref{eq:boundN1}  
and adding the terms to obtain the result of the lemma, namely,
\[
\EE(M_1+M_2)\leq \EE(M_1+M_3+N_2)\leq \frac{1}{p}\EE(M_3+N_2)\leq \frac{1}{2p-1}\left (\frac{2F_{\epsilon}}{h(C)}+ \log_\gamma\left(\frac{C}{\alpha_0}\right)\right).
\]
}

\subsection{Final bound on the expected stopping time}
We finally have the following theorem which trivially follows from Lemmas \ref{lem:hittime1} and \ref{lem:hittime2}. 
%

\begin{theorem} \label{th:mainbound}
Under the condition that $p>1/2$, the hitting time $N_\epsilon$ is bounded in expectation as follows
\[
\EE (N_\epsilon) \leq \frac{2p}{(2p-1)^2}\left (\frac{2F_{\epsilon}}{h(C)}+ \log_\gamma\left(\frac{C}{\alpha_0}\right)\right).
\]
\end{theorem}
\proof{
Clearly 
\[
\EE (N_\epsilon) = \EE\left (\sumn \Lambda_k\right ) + \EE \left (\sumn (1-\Lambda_k)\right )
\]
and, hence, using  Lemmas \ref{lem:hittime1} and \ref{lem:hittime2} we have 
\[
\EE (N_\epsilon)\leq \frac{1}{2p} \EE (N_\epsilon) +   \frac{1}{2p-1} \left (\frac{2F_{\epsilon}}{h(C)}+ \log_\gamma\left(\frac{C}{\alpha_0}\right)\right).
\]
The result of the theorem easily follows. 
}

\paragraph{Summary of our complexity analysis framework.}\label{summary-page}  We have considered a(ny)
algorithm in the framework Algorithm \ref{alg:generic} with probabilistically
sufficiently accurate models as in Definition
\ref{def:random_poised_mart}. 
We have
developed a methodology to obtain (complexity) bounds on the number of
iterations $N_{\epsilon}$ that such an algorithm  takes to reach
desired accuracy. It is important to note that, while we simply provide the bound on $\EE(N_\epsilon)$ it is easy to extend the analysis of the same stochastic processes to provide bounds on $P\{N_\epsilon> K\}$, for any $K$ larger than the bound
 on $\EE(N_\epsilon)$, in particular it can be shown that $P\{N_\epsilon> K\}$ decays exponentially with $K$. 
 
 While in our analysis we assumed that the constant $\gamma$ by which we  decrease and increase $\alpha_k$ is the same, our analysis can be quite easily extended to the case when the constants for increase and decrease are different, say $\gamma_{inc}$ and $\gamma_{dec}$. In this case the threshold on the probability $p$ may no longer be $1/2$ but will be larger if $\gamma_{inc}/\gamma_{dec}<1$ and smaller, otherwise. Some of the constants in the upper bound on  $\EE(N_\epsilon)$ with change accordingly. 
  
Our approach is valid provided  that all of the conditions in Assumption
\ref{ass:alg_behave} hold.  
Next we
show that all these conditions are satisfied by steepest-descent linesearch
methods in the nonconvex, convex and strongly convex case; by general
linesearch methods in the nonconvex case; by cubic regularization
methods (ARC) for nonconvex objectives. In particular, we will specify
what we mean by a probablistically sufficiently accurate first-order
and second-order model in the case of linesearch and cubic
regularization methods, respectively.


\section{The line-search algorithm}\label{sec:ls_algorithm}
We will now apply the generic analysis outlined in the previous
section to the case of the following  simple probabilistic line-search
algorithm.

\begin{algorithm}\label{alg:LS_Random} {\bf A line-search algorithm with random models} \\
\begin{description}
\item[{Initialization}]\ \\
Choose constants $\gamma\in (0,1)$, $\theta\in (0,1)$ and $\alpha_{\max}>0$. Pick initial $x^0$ and $\alpha_0<\alpha_{\max}$. Repeat for $k=0, 1, \ldots$
 \item[{ 1. Compute a model and a step}] \ \\
Compute a random model $m_k$ and use it to generate a direction
$g^k$. \\
Set the step $s^k=-\alpha_k g^k$.

\item[{ 2. Check sufficient decrease}]\ \\  
Check if  
\begin{equation}\label{suff_decrease}
f(x^k-\alpha_k g^k) \leq f(x^k)-\alpha_k\theta\| g^k\|^2.
\end{equation}

\item[{3. Successful step}] \ \\ 
If \eqref{suff_decrease} holds, then $x^{k+1}:=x^k-\alpha_k g^k$ and
$\alpha_{k+1}= \min \{ \alpha_{\max}, \gamma^{-1}\alpha_k\}$. Let $k:=k+1$.

\item[{ 4. Unsuccessful step}]\ \\  
Otherwise, $x^{k+1}:=x^k$, set 
$\alpha_{k+1}=\gamma \alpha_k$. Let $k:=k+1$.
\end{description}
\end{algorithm}

For the linesearch algorithm, the key ingredient is a search direction
selection on each iteration. 
In our case we assume that the search direction is random and
satisfies some accuracy requirement that we discuss below.
The choice of model in this algorithm is a simple  linear model $m_k(x)$, which gives rise to the search direction $g^k$, specifically,
 $m_k(x)=f(x^k)+  (x-x^k)^Tg^k$. We will consider more general
 models in the next section, Section \ref{section:generaldescentmodels}. 


Recall  Definition \ref{def:random_poised_mart}. Here we describe the specific requirement we apply to the models in the case of line search. 

\begin{definition}\label{def:random_good_dir}
We say that a sequence of random models and corresponding directions $\{M_k, G_k\}$ is $(p)$-probabilistically
 "sufficiently accurate" for Algorithm \ref{alg:LS_Random} for a corresponding sequence $\{\Alpha_k, X^k\}$, if there exists a constant $\kappa>0$,  such that  the indicator variables
\[
I_k \; =  \; \mathbbm{1}\{\|G^k-\nabla f(X^k)\|\leq \kappa \Alpha_k \|G^k\|\}
\]
satisfy the following submartingale-like condition
\[
P(I_k=1 | {\cal F}^M_{k-1}) \; \geq \; p,
\]
where ${\cal F}^M_{k-1}=\sigma(M_0,\ldots,M_{k-1})$ is the $\sigma$-algebra generated by $M_0,\ldots,M_{k-1}$.
\end{definition}

As before, each iteration for which $I_k=1$ holds is called a true
iteration. It follows that for every realization of the algorithm,
on  all true iterations, we have 
\begin{equation}\label{eq:grad_approx}
\|g^k-\nabla f(x^k)\|\leq \kappa \alpha_k \|g^k\|,
\end{equation}
which implies, using  $\alpha_k\leq \alpha_{\max}$ and the  triangle
inequality, that
\begin{equation}\label{eq:boundgk}
\|g^k\| \geq \frac {\|\nabla f(x^k)\|}{1+ \kappa \alpha_{\max}}.
\end{equation}

For the remainder of the analysis of Algorithm \ref{alg:LS_Random}, we make the following assumption. 
\begin{assumption} \label{ass:model_qual} 
The sequence of random models and corresponding directions $\{M_k, G_k\}$,  generated in Algorithm \ref{alg:LS_Random}, is $(p)$-probabilistically "sufficiently accurate" for the corresponding random sequence $\{\Alpha_k, X^k\}$, with $p>1/2$.
 \end{assumption}

We also make a standard assumption on the smoothness of $f(x)$ for the remainder of the paper. 
\begin{assumption} \label{ass:Lip_cont} \quad
  $f\in\mathcal{C}^1(\RR^n)$,  is globally bounded below by $f_*$, and has globally Lipschitz continuous gradient $\nabla f$, namely,
\begin{equation}\label{nablaf:Lips}
\|\nabla f(x)-\nabla f(y)\|\leq L\|x-y\|\quad {\rm for\,\, all}\,\, x,\, y\in \RR^n \,\,{\rm and\,\, some}\,\, L>0.
\end{equation}
\end{assumption}

\subsection{The nonconvex case, steepest descent}
As mentioned before, our goal in
the nonconvex case is to compute a bound on the expected number of
iterations $k$ that Algorithm \ref{alg:LS_Random} requires to
obtain an iterate $x^k$ for which $\|\nabla f(x^k)\|\leq \epsilon$. 
We will now compute the specific quantities and expressions defined in
Sections \ref{sec:elem_conv_anal} and \ref{sec:akfkproc}, that allow us
to apply the analysis of our general framework to the specific case of 
Algorithm \ref{alg:LS_Random} for nonconvex functions.

Let $N_\epsilon$ denote, as before,  the number of iterations that are
taken until $\|\nabla f(X^k)\|\leq \epsilon$ occurs (which is a random variable). 
Let us consider the stochastic process $\{\Alpha_k, F_k\}$ with
$F_k=f(x^0)-f(X^k)$ and let $F_{\epsilon}=f(x^0)-f_*$. 
Then $F_k\leq F_{\epsilon}$,  for all $k$.

Next we show that Assumption \ref{ass:alg_behave}  is verified. First we
derive an expression for the constant $C$, related to the size of the
stepsize $\alpha_k$.

\begin{lemma}\label{stepsize_threshold_lemma}
Let Assumption \ref{ass:Lip_cont}  hold. For every realization  of Algorithm \ref{alg:LS_Random},  if iteration $k$ is true (i.e. $I_k=1$), and if 
\begin{equation}\label{alpha_bd}
\alpha_k\leq C= \frac{1-\theta}{0.5L+\kappa},
\end{equation}
then \eqref{suff_decrease} holds. In other words, when \eqref{alpha_bd} holds, any true iteration is also a successful one.
\end{lemma}

\proof{Condition \eqref{nablaf:Lips} implies the following overestimation property for all $x$ and $s$ in $\RR^n$,
\[
f(x+s)\leq f(x)+s^T\nabla f(x) +\frac{L}{2}\|s\|^2,
\]
which implies
\[
\begin{array}{lcl}
\displaystyle f(x^k-\alpha_k g^k)&\leq& f(x^k)-\alpha_k (g^k)^T\nabla f(x^k)+\frac{L}{2}\alpha_k^2\|g^k\|^2.
\end{array}
\]
Applying the Cauchy-Schwarz  inequality and \eqref{eq:grad_approx} we have
\[
\begin{array}{lcl}
\displaystyle f(x^k-\alpha_k g^k)&\leq&  f(x^k)-\alpha_k (g^k)^T[\nabla f(x^k)-g^k]-\alpha_k \|g^k\|^2\left[1-\frac{L}{2}\alpha_k\right]\\

\displaystyle &\leq& f(x^k)+\alpha_k \|g^k\|\cdot\|\nabla f(x^k)-g^k\| -\alpha_k \|g^k\|^2\left[1-\frac{L}{2}\alpha_k\right]\\

\displaystyle &\leq & f(x^k)-\alpha_k \|g^k\|^2\left[1-\left(\kappa+\frac{L}{2}\right)\alpha_k\right].
\end{array}
\]
 It follows that
\eqref{suff_decrease} holds whenever $f(x^k)-\alpha_k \|g^k\|^2[1-(\kappa+0.5L)\alpha_k]\leq f(x^k)-\alpha_k\theta\|g^k\|^2$
which is equivalent to \eqref{alpha_bd}.}

From Lemma \ref{stepsize_threshold_lemma}, and from 
\eqref{suff_decrease} and \eqref{eq:boundgk}, 
for any realization of Algorithm \ref{alg:LS_Random} which gives us
 the specific sequence $\{\alpha_k, f_k\}$, 
the following hold. 
 \begin{itemize}
\item 
If $k$ is a true and successful iteration, then
\[
f_{k+1} \geq f_k + \frac{\theta\|\nabla f(x^k)\|^2\alpha_k}{(1+\kappa\alpha_{\max})^2}
\]
and
\[
\alpha_{k+1}=\gamma^{-1}\alpha_k.
\]

\item
If $\alpha_k\leq C$, where $C$ is defined in \eqref{alpha_bd}, and iteration $k$ is true, then it is also successful. 
\end{itemize}

Hence, Assumption \ref{ass:alg_behave} holds and the process
$\{\Alpha_k, F_k\}$ behaves exactly as our generic process
\eqref{eq:proc1_Zk}-\eqref{eq:proc1_Yk} in Section \ref{sec:akfkproc},
with $C$ defined in \eqref{alpha_bd} and  the specific choice of $h( \Alpha_k)= \frac{\theta\epsilon^2\Alpha_k}{(1+\kappa\alpha_{\max})^2}$.

Finally, we use Theorem \ref{th:mainbound}  and  substituting the expressions for $C$,
  $h(C)$ and $F_{\epsilon}$  into the bound on  $\EE(N_\epsilon)$ we obtain  the following complexity result. 

\begin{theorem}\label{th:finalnoncvxls} 
Let Assumptions \ref{ass:model_qual} and \ref{ass:Lip_cont}  hold. Then
  the expected number of iterations that Algorithm  \ref{alg:LS_Random} takes 
until $\|\nabla f(X^k)\|\leq \epsilon$ occurs is bounded as follows
$$
\EE(N_\epsilon)\leq \frac{2p}{(2p-1)^2}\left [\frac{M}{\epsilon^2}+\log_{\gamma}\left(\frac{1-\theta}{\alpha_0 (0.5L+\kappa)}\right)\right],
$$
where $M=\frac{(f(x^0)-f_*)(1+\kappa \alpha_{\max})^2(0.5L+\kappa)}{\theta(1-\theta)}$ is a constant independent of $p$ and $\epsilon$. 
\end{theorem}

\begin{remark} We note that the dependency of the expected number of iterations on ${\epsilon}$ is of the order  $1/\epsilon^2$, as expected from a line-search method applied to a smooth nonconvex problem. The dependency on $p$ is rather intuitive as well: if $p=1$, then the deterministic complexity is recovered, while as $p$ approaches $1/2$,  the expected number of iterations goes to infinity, since 
the models/directions are arbitrarily bad as often as they are good. 
\end{remark}

\subsection{The nonconvex case, general descent}
\label{section:generaldescentmodels}
In this subsection, we explain how the above analysis of the
 line-search method extends from the nonconvex steepest descent 
case to a general nonconvex descent case. 

In particular, we consider that in Algorithm \ref{alg:LS_Random}, $s^k=\alpha_kd^k$ (instead of $-\alpha_kg^k$),
where $d^k$ is any direction that satisfies the following standard conditions.
\begin{itemize}
\item There exists a constant $\beta>0$, such that 
\begin{equation}\label{gen_d:cond1}
\frac{(d^k)^Tg^k}{\|d^k\|\cdot\|g^k\|}\leq -\beta, \quad \forall k.
\end{equation}
\item There exist  constants $\kappa_1, \kappa_2>0$, such that 
\begin{equation}\label{gen_d:cond2}
\kappa_1\|g^k\|\leq \|d^k\|\leq \kappa_2\|g^k\|, \quad \forall k.
\end{equation}
\end{itemize}

The sufficient decrease condition \eqref{suff_decrease} is replaced by 
\begin{equation}\label{suff_decrease_d}
f(x^k+\alpha_k d^k) \leq f(x^k)+\alpha_k\theta (d^k)^T g^k. 
\end{equation}

It is easy to show that a simple variant of Lemma
\ref{stepsize_threshold_lemma} applies.

\begin{lemma}\label{gen_d:stepsize_threshold_lemma}
Let Assumption \ref{ass:Lip_cont}  hold. Consider Algorithm
\ref{alg:LS_Random} with $s^k=\alpha_k d^k$ and sufficient decrease
condition \eqref{suff_decrease_d}. Assume that $d^k$ satisfies
\eqref{gen_d:cond1} and \eqref{gen_d:cond2}. Then, 
for every realization  of the resulting algorithm,  if iteration $k$ is true (i.e. $I_k$ holds), and if 
\begin{equation}\label{gen_d:alpha_bd}
\alpha_k\leq C=\frac{\beta(1-\theta)}{0.5L\kappa_2+\kappa},
\end{equation}
then \eqref{suff_decrease_d} holds. In other words, when \eqref{gen_d:alpha_bd} holds, any true iteration is also a successful one.
\end{lemma}

\proof{The first displayed equation in the proof of Lemma
  \ref{stepsize_threshold_lemma} provides
\[
\begin{array}{lcl}
\displaystyle f(x^k+\alpha_k d^k)&\leq& f(x^k)+\alpha_k (d^k)^T\nabla f(x^k)+\frac{L}{2}\alpha_k^2\|d^k\|^2.
\end{array}
\]
Applying the Cauchy-Schwarz  inequality, \eqref{eq:grad_approx} and
the conditions  \eqref{gen_d:cond2}  on $d^k$ we have
\[
\begin{array}{lcl}
\displaystyle f(x^k+\alpha_k d^k)&\leq&  f(x^k)+\alpha_k (d^k)^T[\nabla f(x^k)-g^k]+\alpha_k (d^k)^Tg^k
+\frac{L}{2}\alpha^2_k\|d^k\|^2 \\
\displaystyle &\leq& f(x^k)+\alpha_k \|d^k\|\cdot\|\nabla f(x^k)-g^k\| +\alpha_k (d^k)^Tg^k
+\frac{L}{2}\alpha^2_k\|d^k\|^2\\
\displaystyle &\leq& f(x^k)+\alpha^2_k\kappa \|d^k\|\|g^k\|+\alpha_k (d^k)^Tg^k
+\frac{L}{2}\alpha^2_k\kappa_2\|d^k\|\|g^k\|\\
\displaystyle &=& f(x^k) +\alpha_k (d^k)^Tg^k+\alpha_k^2\|d^k\|\|g^k\| \left(\kappa+\kappa_2 \frac{L}{2}\right).
\end{array}
\]
 It follows that
\eqref{suff_decrease_d} holds whenever 
\[
\alpha_k (d^k)^Tg^k+\alpha_k^2\|d^k\|\|g^k\| \left(\kappa+\kappa_2
  \frac{L}{2}\right) \leq \alpha_k\theta (d^k)^Tg^k,
\]
or equivalently, since $\alpha_k>0$, whenever
\[
\alpha_k\|d^k\|\|g^k\| \left(\kappa+\kappa_2
  \frac{L}{2}\right) \leq -(1-\theta) (d^k)^Tg^k.
\]
Using \eqref{gen_d:cond1}, the latter displayed equation holds
whenever $\alpha_k$ satisfies \eqref{gen_d:alpha_bd}.
}

We conclude this extension to general descent directions
 by observing that if  $k$ is a true and successful iteration,  using 
the sufficient decrease condition \eqref{suff_decrease_d}, the
conditions \eqref{gen_d:cond1} and \eqref{gen_d:cond2} on $d^k$ and
\eqref{eq:boundgk}, we obtain that
\[
f_{k+1} \geq f_k + \frac{\theta\kappa_1\beta\|\nabla f(x^k)\|^2\alpha_k}{(1+\kappa\alpha_{\max})^2}.
\]
Hence, Assumption \ref{ass:alg_behave} holds for this case as well and the
remainder of the analysis is exactly the same as for the steepest descent case.

\subsection{The convex case}

We now analyze the expected complexity of Algorithm
\ref{alg:LS_Random} in the case when $f(x)$ is a convex function, that is when the following assumption holds. 

\begin{assumption} \label{ass:convex} \quad $f\in \mathcal{C}^1(\RR^n)$ is convex and has bounded level sets so that
\begin{equation}\label{bounded_levelsets}
\|x-x^*\|\leq D \quad {\rm for\,\,all}\,\,x\,\,{\rm with}\,\,f(x)\leq f(x^0),
\end{equation}
where $x^*$ is a global minimizer of $f$. Let $f^*=f(x^*)$.
\end{assumption}

In this case,  our goal is to bound  the expectation of  $N_\epsilon$ - the  number of iterations taken by Algorithm \ref{alg:LS_Random} until 
\begin{equation}\label{fct_decrease_cond}
f(X^k)-f^*\leq \epsilon
\end{equation}
occurs. We denote $f(X^k)-f^*$ by $\Delta_k^f$ and define  $F_k=\frac{1}{\Delta^f_k}$. Clearly, $N_\epsilon$ is also the number of iterations taken until $F_k  \geq \frac{1}{\epsilon}=F_{\epsilon}$
occurs.

Regarding Assumption \ref{ass:alg_behave}, Lemma
\ref{stepsize_threshold_lemma} provides the value for the
constant $C$, namely, that
whenever $\Alpha_k\leq C$ with 
$C= \frac{1-\theta}{0.5L+\kappa}$,  then every true iteration is also successful.
We  now show that on true and successful
iterations, $F_k$ is increased by at least some function value $h(\Alpha_k)$ for all $k<N_{\epsilon}$. 

\begin{lemma}\label{lem:fctdecrease}
Let Assumptions  \ref{ass:Lip_cont}  and \ref{ass:convex} hold. Consider any realization of Algorithm \ref{alg:LS_Random}. 
For every iteration $k$ that is true and successful, we have
\begin{equation}\label{convex:fctdecrease}
f_{k+1}\geq f_k+ \frac{\theta\alpha_k}{D^2 (1+\kappa\alpha_{\max})^2}.
\end{equation}
\end{lemma}

\proof{Note that convexity of $f$  implies that for all $x$ and $y$,
\[
f(x)-f(y)\geq \nabla f(y)^T(x-y),
\]
and so by using $x=x^*$ and $y=x^k$, we have
\[
-\Delta_k^f=f(x^*)-f(x^k) \geq \nabla f(x^k)^T(x^*-x^k)\geq -D\|\nabla f(x^k)\|,
\]
where to obtain the last inequality, we used Cauchy-Schwarz inequality and \eqref{bounded_levelsets}.
Thus when $k$ is a true iteration, \eqref{eq:boundgk}  further provides
\[
\frac{1}{D}\Delta_k^f\leq \|\nabla f(x^k)\|\leq  (1+\kappa\alpha_{\max})\|g^k\|.
\]
When $k$ is also successful, 
\[
\Delta_k^f-\Delta_{k+1}^f=f(x^k)-f(x^{k+1})\geq \theta\alpha_k\|g^k\|^2\geq \frac{\theta\alpha_k}{D^2 (1+\kappa\alpha_{\max})^2}(\Delta_k^f)^2.
\]
Dividing the above expression by $\Delta_k^f\Delta_{k+1}^f$, we have that on all true and successful iterations 
\[
\frac{1}{\Delta_{k+1}^f}-\frac{1}{\Delta_k^f}\geq \frac{\theta\alpha_k}{D^2 (1+\kappa\alpha_{\max})^2}\frac{\Delta_k^f}{\Delta_{k+1}^f}\geq \frac{\theta\alpha_k}{D^2 (1+\kappa\alpha_{\max})^2},
\]
since $\Delta_k^f\geq \Delta_{k+1}^f$.  Recalling the definition of $f_k$ completes the proof. 
}

Similarly to the nonconvex case, we conclude from  Lemmas
\ref{stepsize_threshold_lemma} and \ref{lem:fctdecrease},
 that for any realization of Algorithm \ref{alg:LS_Random} the following have to happen. 
\begin{itemize}
\item If $k$ is a true and successful iteration, then
\[
f_{k+1} \geq f_k + \frac{\theta\alpha_k}{D^2 (1+\kappa\alpha_{\max})^2}
\]
and
\[
\alpha_{k+1}=\gamma^{-1}\alpha_k.
\]
\item
If $\alpha_k\leq C$, where $C$ is defined in \eqref{alpha_bd}, and iteration $k$ is true, then it is also successful. 
\end{itemize}

Hence, Assumption \ref{ass:alg_behave} holds and the process
$\{\Alpha_k, F_k\}$ behaves exactly as our generic process
\eqref{eq:proc1_Zk}-\eqref{eq:proc1_Yk} in 
Section \ref{sec:akfkproc},  with $C$ defined in \eqref{alpha_bd} and
the specific choice of $h(\Alpha_k)= \frac{\theta\Alpha_k}{D^2
  (1+\kappa\alpha_{\max})^2}$.

Theorem \ref{th:mainbound} can be immediately applied together with  the above expressions for $C$,
$h(C)$ and $F_{\epsilon}$, yielding the following complexity bound.


\begin{theorem}\label{th:finalcvxls} Let Assumptions
  \ref{ass:model_qual}, \ref{ass:Lip_cont} and \ref{ass:convex}
  hold. Then the expected number of iterations that Algorithm \ref{alg:LS_Random} takes 
until $ f(X^k)-f^*\leq \epsilon$ occurs is bounded by 
$$
\EE(N_\epsilon)\leq \frac{2p}{(2p-1)^2}\left [\frac{M}{\epsilon}+\log_{\gamma}\left ( \frac{1-\theta}{\alpha_0 (0.5L+\kappa)}\right )\right],
$$
where $M= \frac{(1+\kappa \alpha_{\max})^2D^2(0.5L+\kappa)}{\theta(1-\theta)}$ is a constant independent of $p$ and $\epsilon$.
\end{theorem}

\begin{remark}
We again note the same dependence on $\epsilon$ in the complexity
bound in Theorem \ref{th:finalcvxls} as in the deterministic convex case and on $p$, as in the nonconvex case. 
\end{remark}

\subsection{The strongly convex case}
We now consider the case of strongly convex objective functions, hence  the following assumption holds.

\begin{assumption}\label{ass:strong_convex} \quad $f\in\mathcal{C}^1(\RR^n)$ is strongly convex, namely, for all $x$ and $y$ and some $\mu>0$,
\[
f(x)\geq f(y)+\nabla f(y)^T(x-y)+\frac{\mu}{2}\|x-y\|^2.
\]
\end{assumption}
Recall our notation $\Delta_k^f=f(X^k)-f^*$. Our goal here is again, as
in the convex case, to bound the expectation on the number of iteration that occur until 
$\Delta_k^f\leq \epsilon$. In the strongly convex case, however, this bound is logarithmic in $\frac{1}{\epsilon}$, just as it is in the case of the deterministic algorithm. 
\begin{lemma}\label{lem:scvx_fct_decrease}
Let Assumption \ref{ass:strong_convex}  hold. Consider any realization of Algorithm \ref{alg:LS_Random}. For every iteration $k$ that is true and successful, we have
\begin{equation}\label{scvx:fct_decrease}
f(x^k)-f(x^{k+1})=\Delta_k^f-\Delta_{k+1}^f\geq \frac{2\mu\theta}{(1+\kappa\alpha_{\max})^2}\alpha_k\Delta_k^f,
\end{equation}
or equivalently,
\begin{equation}\label{scvx:fct_decrease_l}
\Delta_{k+1}^f\leq \left(1-\frac{2\mu\theta}{(1+\kappa\alpha_{\max})^2}\alpha_k\right)\Delta_k^f.
\end{equation}
\end{lemma}

\proof{Assumption \ref{ass:strong_convex} implies, for $x=x^k$ and $y=x^*$, that [see \cite{nesterovtext}, Th 2.1.10]
\[
\Delta_k^f \leq \frac{1}{2\mu}\|\nabla f(x^k)\|^2 
\]
or equivalently,
\[
\sqrt{2\mu \Delta_k^f}\leq \|\nabla f(x^k)\|\leq (1+\kappa\alpha_{\max})\|g^k\|,
\]
where in the second inequality we used \eqref{eq:boundgk}.
The bound \eqref{scvx:fct_decrease} now follows from the sufficient decrease condition  \eqref{suff_decrease}.}

Note that from \eqref{scvx:fct_decrease}  we have that if $\Delta_k^f>0$ and $\alpha_k> (1+\kappa\alpha_{\max})^2/(2\mu\theta)$ 
then the iteration is unsuccessful. Hence, for an iteration to be  successful 
we must have $\alpha_k\leq (1+\kappa\alpha_{\max})^2/(2\mu\theta)$. We also know that a true iteration is successful when
$\alpha_k\leq C$, where $C$ defined in \eqref{alpha_bd}, assuming that $C\leq (1+\kappa\alpha_{\max})^2/(2\mu\theta)$. To simplify the 
analysis we will simply assume that this inequality holds, by an
appropriate choice of the parameters, which can done without loss of generality. 

We now define $F_k=\log\frac{1}{\Delta_k^f}$ and
$F_\epsilon=\log\frac{1}{\epsilon}$, and the  hitting time  $N_\epsilon$ is the number of iterations taken until $\Delta_k^f\leq \epsilon$. 

As in the convex case, using Lemmas \ref{stepsize_threshold_lemma} and
\ref{lem:scvx_fct_decrease}, we conclude that, for any realization of
Algorithm \ref{alg:LS_Random}, the following have to happen. 
\begin{itemize}
\item If $k$ is a true and successful iteration, then
\[
f_{k+1} \geq f_k -\log \left(1-\frac{2\mu\theta}{(1+\kappa\alpha_{\max})^2}\alpha_k\right). 
\]
and
\[
\alpha_{k+1}=\gamma^{-1}\alpha_k.
\]
\item
If $\alpha_k\leq C$, where $C$ defined in \eqref{alpha_bd}, and iteration $k$ is true, then it is also successful. 
\end{itemize}

Hence, again, Assumption \ref{ass:alg_behave} holds  and the process
$\{\Alpha_k, F_k\}$ behaves exactly as our generic process
\eqref{eq:proc1_Zk}-\eqref{eq:proc1_Yk} in Section \ref{sec:akfkproc},
with $C$ defined in \eqref{alpha_bd} and
the specific choice of 
$$h(\Alpha_k)= -\log \left(1-\frac{2\mu\theta}{(1+\kappa\alpha_{\max})^2}\Alpha_k\right).$$


By using the above expressions for $C$, $h(C)$ and $F_{\epsilon}$, again as in the convex case, we have the following complexity bound for the strongly convex case. 

\begin{theorem}\label{th:finalcvxls2} Let Assumptions
  \ref{ass:model_qual},
\ref{ass:Lip_cont}  and \ref{ass:strong_convex}  hold. Then the expected number of iterations that Algorithm
\ref{alg:LS_Random} takes
until $ f(X^k)-f^*\leq \epsilon$ occurs is bounded by 
$$
\EE(N_\epsilon)\leq\frac{2p}{(2p-1)^2}\left [M\log\left(\frac{1}{\epsilon}\right)+\log_{\gamma}\left ( \frac{1-\theta}{\alpha_0 (0.5L+\kappa)}\right )\right], 
$$
where $M= -\log \left(1-\frac{2\mu\theta (1-\theta)}{(1+\kappa\alpha_{\max})^2(0.5L+\kappa)}\right)$ is a constant independent of $p$ and $\epsilon$.
\end{theorem}

 \begin{remark} Again, note the same  dependence of the complexity
   bound in Theorem \ref{th:finalcvxls2} on $\epsilon$
   as for the deterministic line-search algorithm, and the same dependence on
   $p$ as for the other problem classes discussed above.
 \end{remark}

\section{Probabilistic second-order models and cubic regularization methods}\label{sec:arc_algorithm}

In this section we consider a randomized version of second-order methods, whose deterministic counterpart achieves optimal complexity
rate \cite{cagotoPII, cgt38}. As in the line-search case, we show that
in expectation, the same rate of convergence applies as in the
deterministic (cubic regularization) case, augmented by a
term that depends on the probability of having accurate models. Here
we revert back to considering general objective functions that are not necessarily
convex.

\subsection{A cubic regularization algorithm with random models}
Let us now consider a cubic regularization method where the following
model
\begin{equation}\label{cubic}
m_k(x^k+s)= f(x^k)+ s^Tg^k+\frac{1}{2}s^Tb^ks+\frac{\sigma_k}{3}\|s\|^3,
\end{equation}
 is approximately minimized on each iteration $k$ with respect to
$s$,
for some vector $g^k$ and a matrix  $b^k$ and some regularization parameter $\sigma^k>0$.
As before we assume that $g_k$ and $b_k$ are realizations of some random variables $G_k$ and $B_k$, which imply that the model is random and we assume that it is sufficiently accurate with probability at least $p$; the details of this assumption will be given after we state the algorithm. 

The step $s^k$ is computed as in \cite{cagotoPI,cagotoPII} to approximately minimize the
model \eqref{cubic}, namely, it is required to satisfy
\begin{equation}\label{s-calc}
(s^k)^Tg^k+(s^k)^Tb^ks^k+\sigma_k\|s^k\|^3=0\,\,{\rm
  and}\,\,(s^k)^Tb^ks^k+\sigma_k\|s^k\|^3\geq 0
\end{equation}
and
\begin{equation}\label{TCs}
\|\nabla m_k(x^k+s^k)\|\leq \kappa_{\theta} \min\{1,\|s^k\|\} \|g^k\|,
\end{equation}
where $\kappa_{\theta}\in (0,1)$ is a user-chosen constant.

Note that \eqref{s-calc} is satisfied if $s^k$ is the global minimizer
of the model $m_k$ over some subspace; in fact, it is sufficient for
$s^k$ to be the global minimizer of $m_k$ along the line $\alpha
s^k$ \cite{cagotoPII}{\footnote{Note that a recently-proposed cubic regularization variant \cite{bcgmt} can dispense with the approximate global minimization condition altogether while maintaining the optimal complexity bound of ARC. A probabilistic variant of \cite{bcgmt} can be
constructed similarly to probabilistic ARC, and our analysis here can be extended to provide same-order complexity bounds.}}
Condition \eqref{TCs} is a relative termination condition for
the model minimization (say over increasing subspaces) and it is
clearly
satisfied at stationary points of the model; ideally it will be satisfied sooner at least in the early iterations of the algorithm \cite{cagotoPII}.

The probabilistic Adaptive Regularization with Cubics (ARC) framework is presented below.
\begin{algorithm}\label{alg:ARC_Random}{\bf An ARC algorithm with random models }\\

\begin{description}

\item[{Initialization}]\ \\ 
Choose  parameters $\gamma \in (0,1)$, $\theta\in
(0,1)$, $\sigma_{\min}>0$ and $\kappa_{\theta}\in (0,1)$. Pick initial $x^0$ and $\sigma_0>\sigma_{\min}$. Repeat for $k=0, 1, \ldots$,

\item[{ 1. Compute a model }]\ \\  
Compute an approximate gradient  $g^k$ and Hessian $b^k$ and
form the model \eqref{cubic}.

\item[{ 2. Compute the trial step $s^k$}]\ \\
Compute the trial step $s^k$ to satisfy \eqref{s-calc} and \eqref{TCs}.
  
\item[{ 3. Check sufficient  decrease}]\ \\
Compute $f(x^k+s^k)$ and 
\begin{equation}\label{rho}
\rho_k=\frac{f(x^k)-f(x^k+s^k)}{f(x^k)-m_k(x^k+s^k)}.
\end{equation}

\item[{4. Update the iterate}]\ \\
Set 
\begin{equation}\label{newiterate} 
x^{k+1}=\left\{ 
\begin{array}{lr}
x^k+s^k\quad {\rm if}\quad \rho_k\geq \theta &\hfill[k\,\,{\rm successful}]\\
x^k \quad {\rm otherwise} &\hfill[k\,\,{\rm unsuccessful}]
\end{array}
\right.
\end{equation}

\item[{5. Update the regularization parameter $\sigma_k$}]\ \\
Set 
\begin{equation}\label{regularization-update} 
\sigma_{k+1}=\left\{ 
\begin{array}{ll}
\max\{\gamma\sigma_k,\sigma_{\min}\}
\quad {\rm if}\quad \rho_k\geq \theta\\
\frac{1}{\gamma}\sigma_k \quad {\rm otherwise.}
\end{array}
\right.
\end{equation}

\end{description}
\end{algorithm}

\begin{remark}
Typically (see e.g. \cite{cagotoPI}) one would   further refine \eqref{newiterate} and
\eqref{regularization-update} by distinguishing between successful and
very successful iterations, when $\rho_k$ is not just positive but
close to $1$. It is beneficial in the deterministic setting to keep the
regularization parameter unchanged on successful iterations when
$\rho_k$ is greater than $\theta$ but is not close to $1$ and only to decrease it  when $\rho_k$ is substantially larger than $\theta$. 
 For simplicity and uniformity  of our general  framework, 
we simplified the parameter update rule. However, the analysis
presented here can be quite easily extended  to the more general case
by slightly extending the flexibility of the stochastic processes.
In practice it is yet unclear if the same strategy will be beneficial,
as ``accidentally" bad models and the resulting unsuccessful steps may
drive the parameter $\sigma_k$ to be larger than it should be,  and
hence a more aggressive  decrease of $\sigma_k$ may be desired. 
This practical study is a subject of future research. 
\end{remark}

\begin{remark}
We have stated Algorithm \ref{alg:ARC_Random} so that it is as close as
possible to known/deterministic ARC frameworks for ease of reading.
We note however, that it is perfectly coherent with the generic
algorithmic framework, Algorithm \ref{alg:generic}, if one sets
$\alpha_k=1/\sigma_k$ and $\rho_k\geq \theta$ as the sufficient
decrease condition. We will exploit this connection in the analysis
that follows.
\end{remark}

The requirement of sufficient model accuracy considered here is  similar to the definition of  probabilistically fully quadratic models
introduced in \cite{Bandeiraetal2014}, though note that we only
require the second-order condition along the trial step $s^k$.
\begin{definition}\label{def:random_good_mod}
We say that a sequence of random models and corresponding directions $\{M_k\}$ is $(p)$-probabilistically
 "sufficiently accurate" for Algorithm \ref{alg:ARC_Random} if
 there exist constants $\kappa_g$ and $\kappa_H$ such that  for any
 corresponding 
random sequence $\{\Alpha_k=1/\Sigma_k, X^k\}$, the  random indicator variables
\[
I_k \; = \; \mathbbm{1} \{\|\nabla f(X^k)-G^k\|\leq \kappa_g \|S^k\|^2\quad {\rm and}\quad
\|(H(X^k)-B^k)S^k\|\leq \kappa_H \|S^k\|^2\}
\]
satisfy the following submartingale-like condition
\[
P(I_k=1 | F^M_{k-1}) \; \geq \; p,
\]
where $F^M_{k-1}=\sigma(M_0,\ldots,M_{k-1})$ is the $\sigma$-algebra generated by $M_0,\ldots,M_{k-1}$.
\end{definition}


As before, for any realization of Algorithm \ref{alg:ARC_Random},
we refer to iterations $k$ when $I_k$ occurs as {\em true} iterations, and otherwise, as {\it false} iterations.
Hence for all true iterations $k$,
\begin{equation}\label{cr:approxgH} 
\|\nabla f(x^k)-g^k\|\leq \kappa_g \|s^k\|^2\quad {\rm and}\quad
\|(H(x^k)-b^k)s^k\|\leq \kappa_H \|s^k\| ^2.
\end{equation}

For the remainder of the analysis of Algorithm \ref{alg:ARC_Random} we make the following assumption. 
\begin{assumption} \label{ass:model_qual_arc} 
The sequence of random models and corresponding directions $\{M_k, S_k\}$,  generated in Algorithm   \ref{alg:ARC_Random}, is $(p)$-probabilistically "sufficiently accurate" for the corresponding random sequence $\{\Alpha_k=1/\Sigma_k, X^k\}$, with $p>1/2$.
 \end{assumption}

Regarding the possibly nonconvex objective $f$, 
in addition to Assumption \ref{ass:Lip_cont}, we also need the following assumption.

\begin{assumption}\label{ass:Lip_nablaf}\quad $f\in\mathcal{C}^2(\RR^n)$ and has globally Lipschitz continuous Hessian $H$, namely,
\begin{equation}\label{nablaH:Lips}
\|H(x)-H(y)\|\leq L_H\|x-y\|\quad {\rm for\,\, all}\,\, x,\, y\in \RR^n \,\,{\rm and\,\, some}\,\, L_H>0.
\end{equation}
\end{assumption}

\subsection{Global convergence rate analysis, nonconvex case}

The next four lemmas give useful properties of Algorithm
\ref{alg:ARC_Random} that are needed later for our stochastic analysis.
\begin{lemma}[Lemma 3.3 in \cite{cagotoPI}]\label{p-arc:model_decrease}
Consider any realization of Algorithm \ref{alg:ARC_Random}. Then on each iteration $k$
we have
\begin{equation}\label{cr:modeldecrease}
f(x^k)-m_k(x^k+s^k)\geq \frac{1}{6}\sigma_k\|s^k\|^3.
\end{equation}
Thus on every successful iteration $k$, we have
\begin{equation}\label{cr:fctdecrease}
f(x^k)-f(x^{k+1})\geq \frac{\theta}{6}\sigma_k\|s^k\|^3.
\end{equation}
\end{lemma}
\proof{Clearly, \eqref{cr:fctdecrease}
follows from \eqref{cr:modeldecrease} and the sufficient decrease condition \eqref{rho}-\eqref{newiterate}. It
remains to prove \eqref{cr:modeldecrease}.
 Combining the first condition on step $s^k$ in \eqref{s-calc}, with the model expression \eqref{cubic} for $s=s^k$ 
 we can write 
\[
f(x^k)-m_k(x^k+s^k)=\frac{1}{2}(s^k)^TB^ks^k+\frac{2}{3}\sigma_k\|s^k\|^3.
\]
The second condition on $s^k$ in \eqref{s-calc} implies
$(s^k)^TB^ks^k\geq -\sigma_k\|s^k\|^3$, which, when used with the
above  equation, gives \eqref{cr:modeldecrease}.}

\begin{lemma}\label{p-arc:sigmamax}
Let Assumptions \ref{ass:Lip_cont}  and  \ref{ass:Lip_nablaf} hold. 
For any realization of  Algorithm \ref{alg:ARC_Random}, if iteration
$k$ is true (i.e., $I_k=1$),  and if
\begin{equation}\label{sigmamax}
\sigma_k\geq  \sigma_c= \frac{2\kappa_g+\kappa_H +L+L_H}{1-\frac{1}{3}\theta},
\end{equation}
  then iteration $k$ is also successful. 
 
\end{lemma}

\proof{Clearly, if $\rho_k-1\geq 0$, then $k$ is successful by definition. Let us consider 
  the case when $\rho_k<1$; then if $1-\rho_k\leq 1-\theta$, $k$ is
  successful. We have from \eqref{rho}, that 
\[
1-\rho_k=\frac{f(x^k+s^k)-m_k(x^k+s^k)}{f(x^k)-m_k(x^k+s^k)}.
\]
Taylor expansion and triangle inequalities give, for some $\xi^k\in [x^k,x^k+s^k]$,
\[
\begin{array}{l}
f(x^k+s^k)-m_k(x^k+s^k)\\[1ex] 
=[\nabla f(x^k)-g^k]^Ts^k
+\frac{1}{2}(s^k)^T[H(\xi^k)-H(x^k)]s^k +
\frac{1}{2}(s^k)^T[H(x^k)-b^k]s^k-\frac{1}{3}\sigma_k\|s^k\|^3\\[1ex]
\leq \|\nabla f(x^k)-g^k\|\cdot\|s^k\|
+\frac{1}{2} \|H(\xi^k)-H(x^k)\|\cdot \|s^k\|^2 +
\frac{1}{2}\|(H(x^k)-b^k)s^k\|\cdot\|s^k\|-\frac{1}{3}\sigma_k\|s^k\|^3\\[1ex]
\leq \left(\kappa_g+\frac{L_H}{2}+\frac{\kappa_H}{2}
  -\frac{1}{3}\sigma_k\right)\|s^k\|^3=(6\kappa_g + 3L_H+3\kappa_H -2\sigma_k)\frac{1}{6}\|s^k\|^3,
\end{array}
\]
where the last inequality follows from the fact that the iteration is true and hence \eqref{cr:approxgH} holds, and from 
Assumption \ref{ass:Lip_nablaf}. This and \eqref{cr:modeldecrease} now give that
$1-\rho_k\leq 1-\theta$ when $\sigma_k$ satisfies \eqref{sigmamax}.}

Note that for the above lemma to hold $\sigma_c$ does not have to depend on $L$. However, in what follows we will need another condition
on $\sigma_c$, which will involve $L$; hence for simplicity of notation we introduced $\sigma_c$ above to satisfy all necessary bounds.

\begin{lemma}\label{lem:gainontrue}
Let Assumptions \ref{ass:Lip_cont}  and \ref{ass:Lip_nablaf} hold. 
Consider any realization of  Algorithm \ref{alg:ARC_Random}. 
 On each true iteration $k$  we have
\begin{equation}\label{step-length-arc}
\|s^k\|\geq \sqrt{\frac{1-\kappa_{\theta}}{\sigma_k+\kappa_s}\|\nabla f(x^k+s^k)\|},
\end{equation}
where $\kappa_s=2\kappa_g+\kappa_H+L+L_H$.
\end{lemma}
\proof{Triangle inequality, equality $\nabla
  m_k(x^k+s)=g^k+b^ks+\sigma_k\|s\| s$ and condition \eqref{TCs} on $s^k$ together give
\begin{equation}\label{eq:gradbound}
\begin{array}{lcl}
\|\nabla f(x^k+s^k)\|&\leq& \|\nabla f(x^k+s^k)-\nabla m_k(x^k+s^k)\| +
\|\nabla m_k(x^k+s^k)\|\\
 &\leq &\|\nabla f(x^k+s^k)-g^k-b^ks^k\|+\sigma_k\|s^k\|^2 +\kappa_{\theta}\min\{1,\|s^k\|\}\|g^k\|.
\end{array}
\end{equation}
 Recalling Taylor expansion of $\nabla f(x^k)$
\[
\nabla f(x^k+s^k) =\nabla f(x^k)+\int_0^1 H(x^k+ts^k)s^kdt,
\]
and applying triangle inequality, again, we have
\[
\begin{array}{lcl}
\|\nabla f(x^k+s^k)-g^k-b^ks^k\|&\leq& \|\nabla f(x^k)-g^k\|+\\[1ex]
&&\left\|\int_0^1
  [H(x^k+ts^k)-H(x^k)]s^kdt\right\| + \|H(x^k)s^k-b^ks^k\| \\[1ex]
&\leq &\left\{\kappa_g+\frac{1}{2}L_H  +\kappa_H \right\}\|s^k\|^2,
\end{array}
\]
where to get the second inequality, we also used \eqref{cr:approxgH}
and Assumption \ref{ass:Lip_nablaf}.

We can bound $\|g^k\|$ as follows
\[
\|g^k\|\leq \|g^k-\nabla f(x^k)\|+\|\nabla f(x^k)-\nabla f(x^k+s^k)\|
+\|\nabla f(x^k+s^k)\|\leq \kappa_g\|s^k\|^2+L\|s^k\|+\|\nabla f(x^k+s^k)\|.
\]
Thus finally, we can bound all the terms on the right hand side of \eqref{eq:gradbound} in terms of $\|s^k\|^2$ and using the fact that $\kappa_{\theta}\in (0,1)$ we can write 
\[
(1-\kappa_{\theta})\|\nabla f(x^k+s^k)\|\leq (2\kappa_g+\kappa_H+L+L_H+\sigma_k)\|s^k\|^2,
\]
which is equivalent to \eqref{step-length-arc}.
}

\begin{lemma}\label{p-arc:steplength} Let Assumptions
\ref{ass:Lip_cont} and \ref{ass:Lip_nablaf} hold. 
Consider any realization of  Algorithm \ref{alg:ARC_Random}. 
 On each true and successful iteration $k$,  we have
\begin{equation}\label{cr:fctdecrease2}
f(x^k)-f(x^{k+1})\geq  \frac{\kappa_f}{(\max\{\sigma_k,\sigma_{c}\})^{3/2}} \|\nabla f(x^{k+1})\|^{3/2},
\end{equation}
where $\kappa_f:=\frac{\theta}{12\sqrt{2}}(1-\kappa_{\theta})^{3/2}\sigma_{\min}$ and $\sigma_{c}$ is defined in \eqref{sigmamax}.
\end{lemma}
\proof{
Combining Lemma \ref{lem:gainontrue},   inequality  \eqref{cr:fctdecrease} from Lemma \ref{p-arc:model_decrease} and the definition of successful iteration in  Algorithm \ref{alg:ARC_Random}  we have, for all  true and successful iterations $k$,
\begin{equation}\label{cr:fctdecrease4}
f(x^k)-f(x^{k+1})\geq  \frac{\theta}{6}(1-\kappa_{\theta})^{3/2}\frac{\sigma_k}{(\sigma_k+\kappa_s)^{3/2}}\|\nabla f(x^{k+1})\|^{3/2}.
\end{equation}
Using that $\sigma_k\geq \sigma_{\min}$ and that $\kappa_s\leq \sigma_{c}$, \eqref{cr:fctdecrease4} implies \eqref{cr:fctdecrease2}.
}

\paragraph{The stochastic processes and global convergence rate
  analysis}
We are now ready to cast Algorithm \ref{alg:ARC_Random} and its behavior into the
generic stochastic analysis framework of Section
\ref{sec:genralscheme}.

For each realization of Algorithm \ref{alg:ARC_Random}, we define 
$$\alpha_k=\frac{1}{\sigma_k}\quad {\rm  and}\quad f_k=f(x^0)-f(x^k),$$ 
and consider the corresponding  stochastic process
$\{\Alpha_k=1/\Sigma_k, F_k=f(X^0)-f(X^k)\}$.  Let $F_\epsilon
=f(x^0)-f^*$ denote the upper bound on the progress measure $F_k$.

As in the case of the line-search algorithm applied to nonconvex
objectives, we would like to bound the expected number of iterations that
Algorithm \ref{alg:ARC_Random} takes 
until $\|\nabla f(X^k)\|\leq \epsilon$ occurs. Here, however, for technical reasons made clear below, we
count the number of iterations until a successful iteration results in $x^{k+1}$ such that
$\|\nabla f(X^{k+1})\|\leq \epsilon$. Let $N_\epsilon$ denote the (random) index  of such an iteration.
(Clearly,  $N_\epsilon$ thus defined is simply one less than the number of iterations that occur until $\|\nabla f(X^{k+1})\|\leq \epsilon$.)


Regarding Assumption  \ref{ass:alg_behave},   Lemmas \ref{p-arc:sigmamax} 
and \ref{p-arc:steplength} provide that the following must hold for
any realization of Algorithm \ref{alg:ARC_Random}.
\begin{itemize}
\item If $k$ is a true and successful iteration, then
\[
f_{k+1} \geq f_k +  \frac{\kappa_f}{(\max\{\sigma_k,\sigma_{c}\})^{3/2}} \|\nabla f(x^{k+1})\|^{3/2}
\]
and
\[
\alpha_{k+1}=\gamma^{-1}\alpha_k.
\]
\item
If $\alpha_k\leq C=\frac{1}{\sigma_c}$, where $\sigma_c$ is defined in
\eqref{sigmamax},
and iteration $k$ is true, then it is also successful. 
\end{itemize}
Hence, once again, Assumption \ref{ass:alg_behave}  holds and the
process $\{\Alpha_k, F_k\}$ behaves exactly as our generic process
\eqref{eq:proc1_Zk}-\eqref{eq:proc1_Yk} in Section \ref{sec:akfkproc},
with $C=\frac{1}{\sigma_c}=\frac{1-\frac{1}{3}\theta}{2\kappa_g+\kappa_H +L+L_H}$, and the specific choice  
\[
h(\Alpha_k)=\kappa_f(\min\{\Alpha_k,C\})^{3/2} \epsilon^{3/2}.
\]
for all $k<N_\epsilon$.  
Finally, the complexity result again follows from Theorem \ref{th:mainbound} and the expressions for $C$, $h(C)$ and $F_{\epsilon}$.
\begin{theorem}\label{th:finalnoncvxls2} Let Assumptions 
\ref{ass:Lip_cont}, \ref{ass:model_qual_arc}  and
\ref{ass:Lip_nablaf} hold.
Then the expected number of iterations that Algorithm \ref{alg:ARC_Random} takes 
until $\|\nabla f(X^{k+1})\|\leq \epsilon$ occurs is bounded by 
$$
\EE(N_\epsilon)\leq  \frac{2p}{(2p-1)^2}\left (\frac{M}{\epsilon^{3/2}}+
\log\left ( \frac{2\kappa_g+\kappa_H+L+L_H}{\sigma_0(1-1/3\theta)}\right )\right),
$$
where $M=\frac{(f(x^0)-f^*)(2\kappa_g+\kappa_H+L+L_H)^{3/2} }{\kappa_f(1-1/3\theta)^{3/2}}$ is a constant independent of $p$ and $\epsilon$. 
\end{theorem}

\begin{remark} We note that the dependency on $\epsilon$ in the above bound  on the expected number of
  iterations is of the order  $\epsilon^{-3/2}$,
  which is of the same order as for the deterministic
  ARC algorithm and is the optimal rate for nonconvex optimization
  using second order models \cite{cgt38}. 
The dependence on $p$ is, again, the same as in the case of line-search and it  is intuitive.  \end{remark}

\begin{remark} Theorem \ref{thm:liminf} stating that $\liminf k \to \infty  \|\nabla f(X^k)\|=0$ almost surely, holds for Algorithm \ref{alg:ARC_Random} since a similar proof applies. 
\end{remark}

\section{Random models} \label{sec:models}

In this section we will discuss and motivate the definition of probabilistically "sufficiently accurate" models.  In particular, Definition \ref{def:random_good_dir} is a modification of the definition of probabilistically fully-linear models, which is used in \cite{Bandeiraetal2014}.  Similarly, Definition \ref{def:random_good_mod} is similar to that of probabilistically fully-quadratic models in  \cite{Bandeiraetal2014}. These definitions serve to provide properties of the model (with some probability) which are sufficient for first-order (in the case of Definition \ref{def:random_good_dir})
and second-order (in the case of Definition \ref{def:random_good_mod}) convergence rates. 

We will now describe several setting where the models are random and satisfy our definitions. 

\subsection{Stochastic gradients and batch sampling}

In \cite{ByrdChinNocedalWu} an adaptive sample size strategy was proposed in the setting where $\nabla f(x)=\sum_{i=1}^N \nabla f_i(x)$, for large values of $N$.
In this case computing $\nabla f(x)$ accurately can be prohibitive, hence, instead an estimate $\nabla f_S(x)=\sum_{i\in S} \nabla f_i(x)$ is often computed in hopes that it provides a good estimate of the gradient and a  descent direction. It is observed in \cite{ByrdChinNocedalWu} that if sample sets $S_k$ on each iteration ensure   that
\begin{equation}\label{eq:samplegrad}
\|\nabla f_{S_k}(x^k)-\nabla f(x^k)\|\leq \mu \|\nabla f_{S_k}(x^k)\|
\end{equation}
for some $\mu\in (0,1)$, then using a fixed step size
\begin{equation}\label{eq:stepsize}
\alpha_k\equiv \alpha \leq \frac{1-\mu}{L}
\end{equation}
the step $s_k=-\alpha \nabla f_{S_k}(x)$ is always a descent step and the line search algorithm converges with the rate $O(\log(1/\epsilon))$ if $f$ is  strongly convex.
Clearly, condition \eqref{eq:samplegrad} implies that the model $M_k(x)=f(x^k)+\nabla f_{S_k}(x^k)^\top (x-x^k)$ is sufficiently accurate according to Definition \ref{def:random_good_dir} for the given fixed step size $\alpha$. 
Hence Assumption \ref{ass:model_qual} on the models can be viewed as a relaxed version of those in  \cite{ByrdChinNocedalWu}, since we allow the condition \eqref{eq:samplegrad} to fail, as long as 
it fails with probability less than $1/2$, conditioned on the past. Moreover, we analyze the practical version of line search algorithm, with a variable step size, which does not have to remain smaller than $ \frac{1-\mu}{L}$ and we provide convergence rates in convex, strongly convex and nonconvex setting. 

Convergence in expectation of a stochastic algorithm is further shown in \cite{ByrdChinNocedalWu}. In particular, under the assumption  that the variance of $\| \nabla f_i(x)\|$ is bounded for all $i$ and   that $\EE_S[\nabla f_S(x^k)]=\nabla f(x^k)$,  it is shown that, for $X^k$ computed after $k$ steps of stochastic gradient descent with a fixed step size, 
$\EE[f(X^k)]$ converges linearly of $f^*$, when $f(x)$ is strongly convex and if   $|S_k|$ - the size of the sample set $S_k$ - grows exponentially with $k$.

Here, again, our results can be viewed as a generalization of the results in  \cite{ByrdChinNocedalWu}. Indeed, let us assume that $\EE_S[\nabla f_S(x^k)]=\nabla f(x^k)$ for each $x^k$ and let $t_k=|S_k|$ - the size of the sample set $S_k$. Since variance of $\| \nabla f_i(x)\|$ is bounded for all $i$,  we have   that $\EE_{S_k}[\|\nabla f_{{S_k}}(x^k)-\nabla f(x^k)\|]\leq \frac{w}{t_k}$, for some fixed $w$, where the expectation is taken over all random sample sets $S_k$ of size $t_k$.  In other words, the variance of  one sample of the stochastic gradient $\| \nabla f_i(x)\|$ is bounded and 
hence the variance of $\nabla f_{S_k}(x^k)$ decays as the size of $S_k$ increases.  

By Chebychev inequality
\[
Pr\{\|\nabla f_{S_k}(x^k)-\nabla f(x^k)\| >{\min\{1/2,\alpha_k\}\|\nabla f(x^k)\|}\} \leq \frac{ w}{\min\{1/2,\alpha_k\}^2\|\nabla f(x^k)^2\||S_k|}. 
\]
If $\|\nabla f_{S_k}(x^k)-\nabla f(x^k)\|\leq {\min\{1/2,\alpha_k\}\|\nabla f(x^k)\|}$, for a particular $x_k$ and a sample set $S_k$, then by applying triangle inequality we have
\[
\|\nabla f_{S_k}(x^k)-\nabla f(x^k)\|\leq  \frac{ \alpha_k\|\nabla f_{S_k}(x^k)\|}{2}.
\]
Hence the probability of the event $\|\nabla f_S(x^k)-\nabla f(x^k)\|\leq  \frac{\alpha_k\|\nabla f_{S_k}(x^k)\|}{2}$ is at least 
\[
1- \frac{ w}{\min\{1/2,\alpha_k\}^2\|\nabla f(x^k)\|^2|S_k|}\geq 1- \frac{ w}{\min\{1/2,\alpha_k\}^2(1+\alpha_k)^2\|\nabla f_{S_k}(x^k)\|^2|S_k|},
\]
 hence 
as long as $|S_k|$ is chosen sufficiently large, then this probability is greater than $1/2$ and $\nabla f_S(x^k)$ provides us with a probabilistically sufficiently accurate model according to Definition \ref{def:random_good_dir}. Hence the theory described in this paper applies to the case of line search based on stochastic gradient. 
Note that, on top of the results in  \cite{ByrdChinNocedalWu}, we not only analyze line search in nonconvex and convex setting, but also
  show the bound  on the expected number 
of iterations until the desired accuracy is reached, rather than the expected accuracy after a given number of iterations. As we have shown earlier, this  implies 
$\liminf$-type convergence with probability one. 
 Moreover, as mentioned on page \pageref{summary-page}, it is not difficult to extend our analysis to show that the  number 
of iterations until the desired accuracy is reached has exponentially decaying tails. 

Analyzing complexity of methods in this setting in terms of the total number of gradient samples is a subject of some current research \cite{2014pasglyetal}. We leave the exact comparison that can be obtained from our results and those existing in current literature as future research, as this requires defining a sample size selection strategy and possible improvement of our results. 
Similarly, we leave for future research the derivations of the models in this setting that satisfy Definition \ref{def:random_good_mod} for the use within the ARC algorithm.

\subsection{Models based on random sampling of function values}

The motivation behind the notions of probabilistically fully-linear and fully-quadratic models introduced in \cite{Bandeiraetal2014} is based on derivative-free models, which are models based on function values, rather than gradient estimates. We will now show how such models fit into our framework. 

Let us first recall the definition of probabilistically fully-linear and quadratic models and pose it in the terms closest to the ones used in this paper
\begin{definition}\label{def:random_fully_mod}

\begin{enumerate}
\item We say that a sequence of random models  $\{M_k\}$ is $(p)$-probabilistically
fully-linear  if there exists constant $\kappa_g$  such that  for any corresponding random sequence $\Delta_k$, $X^k$, the random indicator variables
\[
I^l_k \; = \mathbbm{1} \; \{\|\nabla f(X^k)-G^k\|\leq \kappa_g \Delta_k \}
\]
satisfy the following submartingale-like condition
\[
P(I^l_k=1 | F^M_{k-1}) \; \geq \; p,
\]
where $F^M_{k-1}=\sigma(M_0,\ldots,M_{k-1})$ is the $\sigma$-algebra generated by $M_0,\ldots,M_{k-1}$.

 \item We call sequence $\{M_k\}$ is $(p)$-probabilistically
fully-quadratic  if there exist constants $\kappa_g$ and $\kappa_H$ such that  for any corresponding random sequence $\Delta_k$, $X^k$, the random indicator variables 
\[
I^q_k \; = \; \mathbbm{1} \{\|\nabla f(X^k)-G^k\|\leq \kappa_g \Delta_k^2\quad {\rm and}\quad
\|H(X^k)-B^k\|\leq \kappa_H \Delta_k\}
\]
satisfy the following submartingale-like condition
\[
P(I^q_k=1 | F^M_{k-1}) \; \geq \; p,
\]
where $F^M_{k-1}=\sigma(M_0,\ldots,M_{k-1})$ is the $\sigma$-algebra generated by $M_0,\ldots,M_{k-1}$.
\end{enumerate}
\end{definition}

The key difference between the conditions in Definition \ref{def:random_fully_mod} and those in  Definitions \ref{def:random_good_dir} and  \ref{def:random_good_mod} is the right hand side of the error bounds  - in
the case of fully-linear and fully-quadratic models $\Delta_k$ is a random variable that does not depend on $M_k$, but in the case of this paper, $\Delta_k$ is replaced by
$\Alpha_k\| G_k\|$ in the case of Definition \ref{def:random_good_dir}  and by $\|S_k\|$ in the case of  Definition \ref{def:random_good_mod}. In other words, the accuracy of the model has to be proportional to the step size which this model produces. Since in  \cite{Bandeiraetal2014} trust region methods are analyzed instead of line search and ARC,   Definition \ref{def:random_fully_mod}  is sufficient. 

Models in \cite{Bandeiraetal2014}  are constructed by sampling
function values in a ball of a given radius around the current iterate
$x^k$ and in all cases  construction of the $k$-th model $M_k$ relies
on the knowledge of the sampling radius.  We will now show that, given
a mechanism of constructing probabilistically fully-linear and fully-quadratic models  for any sequence of
radii (as described in \cite{Bandeiraetal2014}), we can modify our
line search algorithm and ARC algorithm, respectively, and extend the convergence rate analysis to utilize these models. 

\paragraph{Line-search with probabilistically fully-linear models} Let us consider Algorithm \ref{alg:LS_Random}  and corresponding random sequence of iterates $X^k$ and step sizes $\Alpha_k$. If a given model $M_k$ is fully-linear   in $B(X^k, \Alpha_k\Xi_k)$ and $\|G_k\|\geq \kappa_\Delta \Xi_k$, for some positive constant $\kappa_\Delta$,  then model $M_k$ is sufficiently accurate,  according to Definitions \ref{def:random_good_dir}.


To achieve this, for instance, in nonconvex case, for all  $\|\nabla f(X^k)\|\geq \epsilon$ consider $\Xi_k\leq \frac{\epsilon}{2\kappa_g\max\{\Alpha_k, 1\}}$, where $\kappa_g$ is the constant in the definition of fully-linear models. 
Then  any fully-linear model $M_k(x)$ is also sufficiently accurate, simply because 
$\|\nabla f(X^k)-G^k\|\leq \kappa_g \Alpha_k\Xi_k\leq \min \{\Alpha_k, 1\}\frac{\epsilon}{2}$   implies $\|G_k\|\geq  \frac{\epsilon}{2}\geq \Xi_k$. 
Similar bounds can be derived for the convex and strongly convex cases.

Consider the following example of a method that produces probabilistically sufficiently accurate models, based on the arguments above. 
Suppose we are estimating gradients of $f(x)$ by a finite difference scheme using step size $\Alpha_k\Xi_k$, with $\Xi_k$, sufficiently small, and suppose we compute the function values using parallel computations. If some of the computations fail to complete (due to an overloaded processor, say) with some probability and the total probability of having a computational failure  in any of the processors at each iteration is less than $1/2$, conditioned on the past, then we obtain probabilistically sufficiently accurate models. Note that, we do not assume the nature of the computational error, when such error occurs, hence allowing for the gradient estimate to be, occasionally, completely inaccurate.

Another example can be derived from \cite{Bandeiraetal2014}, where  it is shown that sparse gradient and Hessian estimates can be obtained by randomly sampling fewer function values than is needed to construct gradient and/or Hessian by finite differences.  Using this sampling strategy, probabilistically fully-linear and fully-quadratic models can be generated at reduced computation cost. Here again, choosing sampling radius to equal $\Delta_k=\Alpha_k\Xi_k$, with sufficiently small $\Xi_k$  will guarantee that the models are also probabilistically sufficiently accurate. 

We now address a more practical approach, when estimates $\Xi_k$ are not chosen to be small enough a priory, but are dynamically decreased, as another parameter in the algorithm. We will outline how our theory can be extended in this case. 
Consider the following modification of Algorithm \ref{alg:LS_Random}.

\begin{algorithm}\label{alg:LS_Random_TR} {\bf Line-search with probabilistically fully-linear models} \\
\begin{description}
\item[{Initialization}]\ \\
Chose constants $\theta\in (0,1)$, $\gamma \in (0,1)$, $\alpha_{\max}>0$ and $\kappa_\Delta>1$. Pick initial $x^0$ and $\alpha_0<\alpha_{\max}$, $\xi_0$. Repeat for $k=0, 1, \ldots$
 \item[{ 1. Compute a model}] \ \\
Compute a model $m_k$, which is probabilistically fully-linear in $B(x^k, \alpha_k\xi_k)$ and use it to generate a direction $g^k$.
\item[{ 2. Check model accuracy}]\ \\
If $\|g^k\|\geq \kappa_\Delta \xi_k$, then set the step $s^k=-\alpha_k g^k$ and continue to Step 3.\\
Otherwise,  $x^{k+1}=x^k$, $\alpha_{k+1}= \alpha_k$,  $\xi_{k+1}= \xi_k/\kappa_\Delta$, return to Step 1. 
\item[{ 3. Check sufficient decrease}]\ \\  
Check if  
\begin{equation}\label{suff_decrease2}
f(x^k-\alpha_k g^k) \leq f(x^k)-\alpha_k\theta\| g^k\|^2.
\end{equation}

\item[{4. Successful step}] \ \\ 
If \eqref{suff_decrease2} holds, then $x^{k+1}:=x^k-\alpha_k g^k$ and
$\alpha_{k+1}= \min \{\alpha_k/\gamma, \alpha_{\max}\}$.

\item[{ 5. Unsuccessful step}]\ \\  
Otherwise, $x^{k+1}=x^k$. \\
$\alpha_{k+1}=\gamma \alpha_k$. 
\end{description}
\end{algorithm}

In the above algorithm, at each iteration we maintain $\xi_k$, which is expected to be an underestimate of the norm of the descent direction, up to a constant, $\kappa_\Delta$. 
The algorithm then uses $\delta_k=\alpha_k\xi_k$ as the radius for constructing fully-linear models. 
After the model is produced, condition $\|g_k\|\geq\kappa_\Delta\xi_k$ is checked. If this condition holds, the algorithm proceeds 
exactly as the original version, but if this condition fails, then
$\xi_k$ is reduced by a constant ($\kappa_\Delta$ is a practical choice, but any other constant can be used) 
and the iteration is declared to be unsuccessful (hence $x^{k+1}=x^k$), and the step size $\alpha_k$ remains the same. 

Let us consider different possible  outcomes for each iteration $k$ for which $\|\nabla f(x^k)\|\geq \epsilon$. 
From our analysis above, we know that if  $\xi_k\leq \frac{\epsilon}{2\kappa_g\alpha_{\max}}$ and the model $m_k$ is fully linear, 
then $\|g_k\|\geq \xi_k$, hence the model is also sufficiently accurate and the iteration of Algorithm \ref{alg:LS_Random_TR} proceeds as in  Algorithm \ref{alg:LS_Random}.
Since $\xi_k$ is never increased, then, once it is small enough, the analysis of  Algorithm \ref{alg:LS_Random_TR} can be reduced to that of Algorithm \ref{alg:LS_Random}.
Then what remains is to estimate the number of iterations that Algorithm \ref{alg:LS_Random_TR} takes until  $\xi_k\leq \frac{\epsilon}{2\kappa_g\alpha_{\max}}$  or $\|\nabla f(x^k)\|\leq \epsilon$
occurs.

While $\xi_k$ is not sufficiently small, we can have the following outcomes:  1)  $\|g_k\|<\kappa_\Delta \xi_k$, in which case $\xi_k$ is reduced, 2) the model is not fully linear and $\|g_k\|\geq \kappa_\Delta \xi_k$, hence the model may not be sufficiently accurate,  but $\xi_k$ is not reduced and 3) the model is fully-linear  and $\|g_k\|\geq \kappa_\Delta \xi_k$, hence the model
is also sufficiently accurate. Hence with probability at least $p$, $\xi_k$ is reduced or the model is sufficiently accurate. It is possible to extend the definition of our stochastic processes  and their analysis 
to compute the upper bounds on the expected number of iterations
Algorithm \ref{alg:LS_Random_TR} takes until $\|\nabla f(x^k)\|\leq
\epsilon$ occurs. This bound will be increased by adding a constant times the number of iterations it takes to achieve $\xi_k\leq \frac{\epsilon}{2\kappa_g\alpha_{\max}}$, which is $O(\log(1/\epsilon))$. 
Again, similar analysis can be carried out for the cases of convex and strongly convex functions.

{\bf ARC with probabilistically fully-quadratic models} Let us
consider Algorithm \ref{alg:ARC_Random}.
In this case, in the same vein with  line-search, we consider
setting $\Delta_k$ in the Definition \ref{def:random_fully_mod} of
probabilistically fully-quadratic models, to a
sufficiently small value or adjusting it in the run of the algorithm so
as to ensure that when the model is fully-quadratic, it is also
sufficiently accurate (at least asymptotically). We will make these two approaches to the choice of
$\Delta_k$ more precise in what follows. To this end, we need a new
variant of Lemma \ref{lem:gainontrue} for the case of probabilistically fully-quadratic models.

\begin{lemma}\label{lem:gainontrue-fquadratic}
Let Assumptions \ref{ass:Lip_cont}  and \ref{ass:Lip_nablaf} hold. 
Consider any realization of  Algorithm \ref{alg:ARC_Random} where we
generate models that are p-probabilistically fully-quadratic according
to Definition \ref{def:random_fully_mod}. Then on each iteration $k$ in
which $I_k^q=1$, we have 
\begin{equation}\label{step-length-arc-fquadratic}
(1-\kappa_{\theta})\|\nabla f(x^k+s^k)\|\leq (2\kappa_g+\kappa_H)
\delta_k\max\{\delta_k,1\}+ (L+L_H+\sigma_k)\|s^k\|^2.
\end{equation}
In particular, if $\epsilon\in (0,1]$, $\max\{L,L_H\}\geq 1$, and
\begin{equation}\label{ARC-fquad-delta}
\delta_k\leq \frac{(1-\kappa_{\theta})\epsilon}{\max\left\{ 2(2\kappa_g+\kappa_H), L+L_H+\sigma_k \right\}},
\end{equation}
then on each iteration $k$ with $\|\nabla f(x^k+s^k)\|\geq \epsilon$ and in
which $I_k^q=1$, we have $\|s^k\|\geq \delta_k$. 

Assume now that $\delta_k=\frac{\xi_k}{\sigma_k}$. Then if
$\epsilon\in (0,1]$, $\max\{L,L_H\}\geq 1$, and
\begin{equation}\label{ARC-fquad-delta2}
\xi_k\leq \frac{(1-\kappa_{\theta})\sigma_{\min}\epsilon}{\max\left\{ 2(2\kappa_g+\kappa_H), L+L_H+\sigma_{\min} \right\}}:=\xi_\epsilon,
\end{equation}
then on each iteration $k$ with $\|\nabla f(x^k+s^k)\|\geq \epsilon$ and in
which $I_k^q=1$, we have $\|s^k\|\geq \delta_k$. 
\end{lemma}
\proof{It follows from Definition \ref{def:random_fully_mod} that 
on each realization of Algorithm \ref{alg:ARC_Random}, we have
\begin{equation}\label{IQ}
\|\nabla f(x^k)-g^k\|\leq \kappa_g \delta_k^2\quad {\rm and}\quad
\|H(x^k)-b^k\|\leq \kappa_H \delta_k
\end{equation}
The proof of \eqref{step-length-arc-fquadratic} now follows
  identically to the proof of Lemma \ref{lem:gainontrue} if one uses
  \eqref{IQ} instead of \eqref{cr:approxgH}.

The choice of $\delta_k$ in \eqref{ARC-fquad-delta} implies
$\delta_k\leq 1$ and so $\|s^k\|\geq \delta_k$ trivially holds when
$\|s^k\|\geq 1$. When  $\|s^k\|<1$, $\|\nabla f(x^k+s^k)\|\geq
\epsilon$, and $\delta_k\leq 1$, \eqref{step-length-arc-fquadratic}
implies 
\[
(1-\kappa_{\theta})\epsilon- (2\kappa_g+\kappa_H)
\delta_k\leq (L+L_H+\sigma_k)\|s^k\|.
\]
Now the condition \eqref{ARC-fquad-delta} on $\delta_k$ implies
$(L+L_H+\sigma_k)\|s^k\|\geq
(1-\kappa_{\theta})\epsilon/[2(2\kappa_g+\kappa_H)]$. Applying again
the upper bound on $\delta_k$ provides $\|s^k\|\geq \delta_k$.

Finally, if $\delta_k=\frac{\xi_k}{\sigma_k}$, and using $\sigma^k\geq
\sigma_{\min}$ for all $k$ due to the algorithm construction,
\eqref{ARC-fquad-delta2} implies \eqref{ARC-fquad-delta}.}

The second part of Lemma \ref{lem:gainontrue-fquadratic} provides that if
p-probabilistically fully-quadratic models are generated with 
$\delta_k$ chosen sufficiently small so that \eqref{ARC-fquad-delta}
holds, then the models are also p-probabilistically sufficiently
accurate. Thus Algorithm \ref{alg:ARC_Random} can be run with 
models sampled in this way and the analysis carries through as
before. For example, as in the case of linesearch, $g^k$ and $b^k$
could be generated by (sufficiently accurate)  finite-difference
schemes using function values, where computations are
done in parallel and where the total probability of computational
failure in any of the processors at each iteration is less than $1/2$.

Note however, that the bound that dictates the choice of a suitably
small $\delta_k$ depends on problem constants that may not be known a
priori. Thus it would be better -- and computationally more efficient
-- to adjust $\delta_k$ during the run of Algorithm \ref{alg:ARC_Random}. 
A modification of Algorithm \ref{alg:ARC_Random} that allows this is
given next, and can be viewed as the analogue for ARC of the line-search
Algorithm \ref{alg:LS_Random_TR}.

\begin{algorithm}\label{alg:LS_Random_ARC} {\bf  ARC with
    probabilistically fully-quadratic models} 
\begin{description}
\item[{Initialization}]\ \\
Choose  parameters $\sigma_{\min}>0$, $\gamma\in (0,1)$, $\theta\in (0,1)$, $0<\kappa_{\theta}<1$ and $\kappa_\Delta>1$. Pick a starting
point $x^0$, a starting value $\sigma_0>\sigma_{\min}$ and $\xi_0>0$. Repeat for $k=0, 1, \ldots$,

\item[{1. Compute a model }] \ \\
Compute a model which is probabilistically fully-quadratic in
$B\left(x^k, \frac{\xi_k}{\sigma_k}\right)$, 
and hence generate  approximate gradient  $g^k$ and Hessian $b^k$.

\item[{2. Compute the trial step $s^k$}]\ \\
Compute the trial step $s^k$ to satisfy \eqref{s-calc} and \eqref{TCs}.

\item[{3. Check model accuracy }] \ \\
If  $\|s^k\|\geq \kappa_{\Delta}\delta_k:=\xi_k/\sigma_k$, then go to
Step 4.\\
Otherwise, set $x^{k+1}=x^k$, $\sigma_{k+1}=\sigma_k$,
$\xi_{k+1}=\xi_k/\kappa_{\Delta}$, and return to Step 1.  

\item[{4. Check sufficient  decrease}]\ \\
Compute $f(x^k+s^k)$ and 
\[
\rho_k=\frac{f(x^k)-f(x^k+s^k)}{f(x^k)-m_k(x^k+s^k)}.
\]

\item[{5. Update the iterate}] \ \\
Set 
\[
x^{k+1}=\left\{ 
\begin{array}{lr}
x^k+s^k\quad {\rm if}\quad \rho_k\geq \theta &\hfill[k\,\,{\rm successful}]\\
x^k \quad {\rm otherwise} &\hfill[k\,\,{\rm unsuccessful}]
\end{array}
\right.
\]

\item[{6. Update the regularization parameter $\sigma_k$}] \ \\
Set 
\[
\sigma_{k+1}=\left\{ 
\begin{array}{ll}
\max\{\gamma\sigma_k,\sigma_{\min}\}
\quad {\rm if}\quad \rho_k\geq \theta\\
\frac{1}{\gamma}\sigma_k \quad {\rm otherwise.}
\end{array}
\right.
\]

\end{description}
\end{algorithm}

Algorithm \ref{alg:LS_Random_ARC} updates $\xi_k$ in order to obtain
an underestimate $\delta_k:=\xi_k/\sigma_k$  on the length of the step $s^k$.
It constructs probabilistically fully-quadratic models in
$B(x^k,\xi_k/\sigma_k)$ and checks whether $\|s^k\|\geq
\kappa_{\Delta}\delta_k$. If that is the case, then the iteration of
the above algorithm
proceeds as (Algorithm \ref{alg:ARC_Random}) before;  note that then, if the model is
fully quadratic then it is also sufficiently accurate. Otherwise, if
the step is too short, then $\xi_k$ is decreased by $\kappa_{\Delta}$,
$x^k$ and $\sigma^k$ remain unchanged and a new model is generated (within
the smaller ball).

Let us consider the behavior of Algorithm \ref{alg:LS_Random_ARC}
while $\|\nabla f(x^k+s^k)\|\geq \epsilon$. It follows from the last
part of Lemma \ref{lem:gainontrue-fquadratic} that since $\xi_\epsilon$ is
independent of $k$ and $\xi_k$ is never increased in the algorithm,  if
$\kappa_\Delta\xi_j\leq \xi_\epsilon$ for some $j$, then $\xi_k$ will remain below
this threshold for all subsequent iterations $k\geq j$; from this $j$
onwards, whenever the model is fully quadratic,  then $\|s^k\|\geq
\kappa_\Delta\delta_k$ and the model is also sufficiently
accurate. Thus from iteration $j$ onwards, Algorithm
\ref{alg:LS_Random_ARC} reduces to Algorithm \ref{alg:ARC_Random} and
the complexity analysis is the same as before. It remains to estimate
the size of $j$, namely, the number of iterations Algorithm
\ref{alg:LS_Random_ARC} takes until $\xi_k\leq \xi_\epsilon$ or
$\|\nabla f(x^k+s^k)\|<\epsilon$. 

Similarly to the linesearch analysis of possible outcomes above, 
we can argue that while $\xi_k$ is not  sufficiently small, at
least with probability $p$, $\xi_k$ is reduced or the model is
sufficiently accurate. Thus, extending our earlier ARC analysis (and
definitions of stochastic processes, etc) to account for the $\xi_k$
updates as well, we would find that the complexity bound for Algorithm
\ref{alg:LS_Random_ARC} is essentially that of Algorithm \ref{alg:ARC_Random} 
plus a $\mathcal{O}(\log(1/\epsilon))$ term (coming from $\log(\xi_0/\xi_\epsilon)\log\kappa_\Delta$) that accounts for the
number of iterations to drive $\xi_k$ below $\xi_\epsilon$.


\section{Conclusions}

We have proposed a general algorithmic framework with random models
and a methodology for analyzing its complexity that relies on bounding
the hitting time of a nondecreasing stochastic process that measures
progress towards optimality. Our framework accounts for linesearch and
cubic regularization methods, for example, and we particularize our results to
obtain precise complexity bounds in the case of nonconvex and convex
functions. Despite allowing our models to be arbitrarily inaccurate sometimes, the
bounds we obtained match their deterministic counterparts in the order
of the accuracy $\epsilon$. The effect of model inaccuracy is reflected by 
the constant multiple of the bound, which is a function of the probability that the
model is sufficiently accurate. We have also briefly discussed ways to
obtain probabilistically sufficiently accurate models as required by
our framework.

The results in the paper assume that the objective $f$ is
deterministic. Obtaining global rates of convergence results for
similar algorithmic frameworks when $f$ is stochastic is a topic of future research. Also,
further exploring ways to efficiently generate probabilistically sufficiently
accurate models may increase the applicability of our results to a
diverse set of problems.

\paragraph*{Acknowledgements}

We would like to thank Alexander Stolyar for helpful discussions on stochastic processes. We also would like to thank Zaikun Zhang, who was instrumental in helping us significantly simplify  the analysis of the stochastic process in Section \ref{sec:genralscheme}.

\bibliographystyle{siam}

\bibliography{ref-random_ls}

\footnotesize


\end{document}